\documentclass[hidelinks,onefignum,onetabnum,onethmnum,onealgnum,oneeqnum]{siamart220329}
\usepackage{graphicx}
\usepackage{subcaption}
\usepackage{gensymb}
\usepackage{algorithm}

\usepackage{lipsum}
\usepackage{amsfonts}
\usepackage{amsmath}
\usepackage{amssymb}
\usepackage{txfonts}
\usepackage{graphicx}
\usepackage{epstopdf}
\usepackage{xcolor}
\newcommand{\red}[1]{{\iffalse #1 \fi}}

\usepackage{algorithmic}
\ifpdf
  \DeclareGraphicsExtensions{.eps,.pdf,.png,.jpg}
\else
  \DeclareGraphicsExtensions{.eps}
\fi

\newcommand{\papertitle}{Riemannian Optimization on Tree Tensor Networks with Application in Machine Learning}

\newsiamremark{remark}{Remark}
\newsiamremark{hypothesis}{Hypothesis}
\crefname{hypothesis}{Hypothesis}{Hypotheses}
\newsiamthm{claim}{Claim}

\headers{Riemannian Optimization on TTNs}{M. Willner, M. Trenti, and D. Lebiedz}

\title{\papertitle}

\newcommand{\TAIS}{Tensor AI Solutions GmbH, Pfaffenhofen a.d. Roth, Germany}

\newcommand{\Ulm}{Institute for Numerical Mathematics, University of Ulm, Ulm, Germany}

\newcommand{\Augs}{Chair of Mathematical Data Science, University of Augsburg, Augsburg, Germany}

\usepackage{hyperref}
\usepackage{cleveref}

\author{
Marius Willner 
\footnote{\label{Augs}\Augs} 
\footnote{\label{TAIS}\TAIS} \and 
Marco Trenti 
\footref{TAIS} \and 
Dirk Lebiedz 
\footnote{\label{Ulm}\Ulm}
} 

\usepackage{amsopn}

\definecolor{mygreen}{rgb}{0,0.5,0}

\usepackage{nicefrac}
\usepackage{dsfont}
\usepackage{amsfonts}
\usepackage{amsmath}
\usepackage{tikz-cd}
\usepackage{amssymb}
\usepackage{txfonts}

\newcommand{\dd}{\mathrm{d}}

\newcommand{\cartesian}{\equiv}
\DeclareMathAlphabet{\mathcal}{OMS}{cmsy}{m}{n}

\newcommand{\N}{\mathbb{N}}

\newcommand{\R}{\mathbb{R}}

\newcommand{\Ortho}{\mathrm{O}}
\newcommand{\Skew}{\mathrm{Skew}}
\newcommand{\St}{\mathrm{St}}

\newcommand{\Tr}[1]{\mathrm{Tr}[#1]}
\newcommand{\HO}[1][x]{H_{#1}^\times\mathcal{M}}
\newcommand{\HM}[1][x]{H_{#1}\mathcal{M}}
\newcommand{\HSt}[1][x]{H_{#1}^{\mathbf{\cartesian}}\mathcal{T}}
\newcommand{\TM}[1][x]{T_{#1}\mathcal{M}}
\newcommand{\VM}[1][x]{V_{#1}\mathcal{M}}
\newcommand{\TMG}[1][x]{T_{[#1]}\mathcal{M/G} }
\newcommand{\TT}[1][x]{T_{#1}\mathcal{T}}
\newcommand{\VT}[1][x]{V_{#1}\mathcal{T}}
\newcommand{\HT}[1][x]{H_{#1}^\times\mathcal{T}}
\newcommand{\TTG}[1][x]{T_{[#1]}\mathcal{T/G} }
\newcommand{\mskew}{\mathrm{skew}}
\newcommand{\grad}{\mathrm{grad}}
\newcommand{\Hess}{\mathrm{Hess}}
\newcommand{\Ret}{\mathrm{R}}

\newcommand{\dotpro}[2]{\langle \, #1 \, | \, #2 \, \rangle}
\newcommand{\scapro}[3]{\langle\, #1 \,|\, #2 \,|\, #3 \,\rangle}

\newcommand{\vect}[1]{\mathrm{vec}(#1)}
\newcommand{\mat}[2]{\mathbf{#1}^{(#2)}}
\newcommand{\lift}{\mathrm{lift}}
\newcommand{\Proj}{\mathrm{Proj}}

\newcommand{\otimesk}{\otimes_\mathrm{K}}
\newcommand{\tl}{{t_\mathrm{L}}}
\newcommand{\tr}{{t_\mathrm{R}}}
\newcommand{\tb}{t_\mathrm{B}}
\newcommand{\T}{\mathrm{T}}

\ifpdf
\hypersetup{
  pdftitle={},
  pdfauthor={}
}
\fi

\usepackage{tikz}
\usetikzlibrary{arrows, matrix, fit, backgrounds, positioning}

\begin{document}

\maketitle

\begin{abstract}
Tree tensor networks (TTNs) are widely used in low-rank approximation and quantum many-body simulation.
In this work, we present a formal analysis of the quotient geometry underlying the TTN parameter space.
Our framework allows for arbitrary horizontal distributions,
and we develop efficient first- and second-order optimization algorithms 
that exploit this geometry. Additionally, we devise a backpropagation
algorithm for training TTNs in a kernel learning setting.
We validate our methods through numerical experiments on a 
representative digit classification task and reveal an important tradeoff
between two different horizontal distributions that are available for TTNs: 
while one offers cleaner geometric statements, the other ultimately leads to more 
efficient algorithms.

\end{abstract}



\section{Introduction}
Tensor networks, also known as low-rank tensor formats, have long been used as a model 
order reduction technique to achieve low-rank approximations of high-dimensional tensors, see e.g. \cite{Grasedyck:2013, Nouy:2017,Hackbusch:2019, Bachmayr:2023}
for overviews. 
They are also popular in quantum many-body simulations \cite{White:1992, Vidal:2003, McCulloch:2007}, where the density matrix renormalization group (DMRG) algorithm is widely used 
for optimization tasks, such as finding ground states of quantum systems. 

More recently, tensor networks have been adapted to data modeling tasks \cite{Steinlechner:2014, DaSilva:2015, Grelier:2019} 
 and specifically to non-linear kernel learning on tensor trains \cite{Stoudenmire:2016, Cheng:2019, Klus:2019, Novikov:2017, Kirstein:2022}, 
with modifications to the DMRG algorithm to accommodate these new applications. 
This approach involves embedding input vectors in a feature space
and using the tensor network to represent a trainable high-dimensional weight tensor.
Both tensor trains (TT) and tree tensor networks (TTN) can be used for the task, but in this work, we will 
focus on the latter, as they are seen to better capture
long-range and strongly correlated systems \cite{Murg:2010,Seitz:2023, Hao:2024, Ceruti:2025}.
In previous work on the topic,
block-coordinate descent schemes were utilized to train TTNs \cite{Cheng:2019,Liu:2019,Hao:2024}, 
leaving full descent methods undiscussed.

Furthermore, it is well known that parameter spaces of TTs and TTNs give rise to quotient manifolds 
\cite{Holtz:2010,Uschmajew:2013, Haegeman:2014}, which are used in particular to
show that low-rank TT and TTN tensors form smooth manifolds embedded in
high-dimensional tensor spaces. Directly working with the geometry of those low-rank manifolds is a popular approach 
in a wide range of applications \cite{Lubich:2013, Steinlechner:2016, Uschmajew:2020, Bachmayr:2023}, 
as opposed to optimizing on the quotient spaces themselves~\cite{DaSilva:2015,Hauru:2021}.

Quotient optimization crucially requires a horizontal distribution on the parameter space, which allows one to numerically
grasp operators and operands of the quotient. The horizontal distribution
is canonically chosen orthogonal to the fibers accompanying the quotient, because this leads to elegant geometric
statements, both for first- and second-order algorithms \cite{Lee:2019, Boumal:2023, Gallier:2020}. 
As opposed to this \textit{orthogonal horizontal space} $H^\times$, Riemannian optimization on 
low-rank tensor manifolds has brought forth the non-orthogonal \emph{Cartesian horizontal space} $H^\equiv$ which forms as a Cartesian product of Stiefel horizontal spaces.
One of its main purposes is to achieve a numerically tangible representation of tangent vectors in high-dimensional tensor spaces \cite{Holtz:2010, Steinlechner:2016}, 
and initial efforts in developing first-order tools for quotient optimization involving this particular horizontal distribution have been reported in \cite{DaSilva:2015}.

\textbf{Summary of contributions \& outline.}
In this work, we develop first- and second-order optimization algorithms on the quotient manifold of tree tensor network parameters and apply them to a machine learning task.
We summarize our contributions as follows:
\begin{itemize} \setlength{\itemsep}{0pt}
    \item We introduce a quotient framework that is versatile enough to capture both horizontal distributions $H^\times$ and $H^\equiv$, 
thereby bridging the gap between the various notions 
used across different communities and enabling direct comparison among them. We address an open question raised in \cite{DaSilva:2015} by constructing the orthogonal horizontal space $H^\times$ of TTNs, 
    as well as projectors onto it.
    \item By allowing for quotient connections with torsion we develop novel second-order tools, 
    that are compatible with arbitrary horizontal spaces and in particular the Cartesian horizontal space. 
    In this vein,
    we generalize results established for optimization on matrix manifolds \cite{Absil:2008,Wen:2010, Absil:2012, Boumal:2023},
    which mainly cover tools related to the orthogonal horizontal space.
    \item Employing the kernel-learning approach of \cite{Stoudenmire:2016, Novikov:2017}, 
        we formulate a backpropagation algorithm on tree tensor networks, 
        which serves as a cornerstone for our machine learning application. 

    \item We provide a comprehensive numerical comparison of the two different horizontal distributions that are available for TTNs. 
While optimization involving the orthogonal horizontal space $H^\times$ better respects the Riemannian submersion structure and 
offers stronger convergence results, we observe no disadvantage in terms of convergence behaviour
for the Cartesian horizontal space $H^\equiv$, which ultimately leads to more efficient algorithms. The difference in computation time is especially pronounced for 
second-order methods, where algorithms using the Cartesian horizontal space $H^\equiv$ benefit from cheaper projections and 
can sometimes skip update steps entirely for large parts of the tree.
\end{itemize}

In \Cref{sec:manifold}, we present our quotient framework, 
two concrete horizontal distributions for TTNs and projectors adhering to them. 
\Cref{sec:opt_tools} reports on first- and second-order quotient optimization tools for TTNs that apply
to arbitrary horizontal spaces.
We also present some efficient retractions available on tree tensor network manifolds, which are cheaper and more easily parallelized
than e.g. the HOSVD retraction \cite{Steinlechner:2014}. 
In \Cref{sec:opt_algs}, the developed tools are employed to form concrete and well-known optimization algorithms, 
namely Riemannian gradient descent, Newton's method and trust-region \cite{Absil:2008, Boumal:2023}. 
We modify these accordingly 
to suit the geometry
adhering to tree tensor networks. We qualitatively compare quotient optimization approaches to optimizing 
directly on the low-rank manifolds and shed light on a peculiarity involving the Cartesian horizontal space: due to the gauge-freedom 
of TTNs, some nodes need not be touched at all during optimization.
\Cref{sec:ml_ttn} covers the theoretical and practical application of quotient optimization on TTNs
to machine learning.
Finally, \Cref{sec:results} showcases a numerical comparison of the two different horizontal distributions $H^\times$ and $H^\equiv$, 
when applied to a digit classification task.  Additionally, we confirm competitiveness of TTN machine learning models 
with other classification
models in terms of generalization power.

\section{Preliminaries}
\label{sec:preliminaries} In this section, we summarize some well-known concepts from multilinear algebra, 
and provide a brief review of tree tensor networks. 
We denote vectors by lower case letters (e.g. $u, v, w$), matrices by upper case letters (e.g. $A, B, C$) 
and tensors by bold upper case letters (e.g. $\textbf{A}, \textbf{B}, \textbf{C})$. 
The identity matrix is denoted by $I_n \in \R^{n \times n}$ and the zero matrix is denoted by $0_n \in \R^{n \times n}$. 
Subscripts are dropped if the dimension is clear.

\subsection{Tensors \& operations} \label{sec:tensors}
Consider real vector spaces $\R^{n_j}$, $n_j \in \N$, $j \in D = \{1, \ldots, d\}$. 
For the standard bases $\{ e^{j}_i, 1 \leq i \leq n_j \}$ of $\R^{n_j}$, 
$\mathcal{X} \in \R^{n_1} \otimes \ldots \otimes  \R^{n_d}$ admits the basis representation
\begin{equation*} 
    \mathcal{X} = \sum_{i_1 = 1}^{n_1} \cdots \sum_{i_d = 1}^{n_d} \mathbf{X}_{i_1\ldots i_d} e^{1}_{i_1} \otimes \ldots \otimes e^{d}_{i_d}
\end{equation*}
with the scalar components $\mathbf{X}_{i_1\ldots i_d}$ forming the entries of a $d$-dimensional array $\mathbf{X} \in \R^{n_1 \times \ldots \times n_d}$.
Throughout this paper, we will mainly work with elements in $\R^{n_1 \times \ldots \times n_d} \cong \R^{n_1} \otimes \ldots \otimes  \R^{n_d}$. 

We denote the (col-major) vectorization of $\mathbf{X} \in \R^{n_1 \times \ldots \times n_d}$ as $\vect{\mathbf{X}} \in \R^{n_1 n_2 \ldots n_d}$. 
The Kronecker product is symbolized by $\otimesk$, in order to
distinguish it from the tensor product $\otimes$. We adopt a modified ordering of Kronecker factors, such that $a \otimesk b = \vect{a \otimes b}$ for 
two vectors $a, b\in \R^n$.
With an ordered subset of dimensions $t = (t_1, \ldots, t_k) \subset D$ and $s = (s_1, \ldots, s_{d-k})$ the increasing ordering of $D \setminus t$, 
we denote the (col-major) $t$-matricization of $\mathbf{X}$ as 
$\mat{X}{t} \in \R^{n_{t_1}\ldots{n_{t_k}\times n_{s_1}\ldots n_{s_{d-k}}}}$.
When working with order-three tensors, 
we will also frequently use non-bold symbols to denote $(1,2)$-matricizations of 
the corresponding tensor, e.g. 
$B := \mat{B}{1,2}$.
For tensors $\mathbf{X}, \mathbf{Y} \in \R^{n_1 \times \ldots \times n_d}$, 
define their inner product as $\dotpro{\mathbf{X}}{\mathbf{Y}} := \vect{\mathbf{X}}^\T\vect{\mathbf{Y}}$. 
The inner product is invariant under matricizations: $\dotpro{\mathbf{X}}{\mathbf{Y}}  = \dotpro{\mat{X}{t}}{\mat{Y}{t}}$. 
Given a ($d-1$)-element ordering $t = (t_1, \ldots, t_{d-1})$ of $D \setminus \{i\}$, define their $t$-contraction \cite[Def. 2]{DaSilva:2015} as
\begin{equation*}
    \dotpro{\mathbf{X}}{\mathbf{Y}}_t = (\mat{X}{t})^\T\mat{Y}{t} \in \R^{n_i \times n_i}.
\end{equation*}
Given a matrix $A_i \in \R^{m_i \times n_i}$, we denote the mode-$i$ product with $\mathbf X$ by $A_i \times_i \mathbf{X}$. 
For multi-mode multiplications with a tuple of matrices $(A_i)_{i \in D}, A_i \in \R^{m_i \times n_i}$, we employ the notation
$$
    (A_1 \otimes \ldots \otimes A_d) \mathbf X = A_1 \times_1 \cdots \, A_d \times_d \mathbf X,
$$
which is understood to be consistent with Kronecker and matricization ordering:
\begin{equation*}
    [(A_1 \otimes \ldots \otimes A_d) \mathbf X]^{(t)} = (A_{t_1} \otimes_K \ldots \otimes_K A_{t_k}) \mathbf X^{(t)} (A_{s_1} \otimes_K \ldots \otimes_K A_{s_{d-k}})^\T.
\end{equation*}
In combination with inner products or tensor contractions with tensor $\mathbf{Y} \in \R^{m_1 \times \ldots \times m_d}$ we write
\begin{equation*}
    \dotpro{\mathbf{Y}}{(A_1 \otimes \ldots \otimes A_d)\mathbf{X}}
    = \scapro{\mathbf{Y}}{ A_1 \otimes \ldots \otimes A_d }{ \mathbf{X}}
    = \dotpro{(A_1^\T \otimes \ldots \otimes A_d^\T)\mathbf{Y}}{\mathbf{X}}.
\end{equation*}

\subsection{Tree tensor networks}
In this subsection, we provide a quick recapitulation of the hierarchical Tucker format (HT) \cite{Hackbusch:2009}. 
Our notation goes along the lines of \cite[Sec. 3]{Uschmajew:2013}, 
however we only allow for order-three tensors in the tree, 
as it better fits descriptions of 
tree tensor networks given in the physics community \cite{Shi:2006, Cheng:2019, Hao:2024}
and simplifies notation later on.

\emph{Structure \& labels.} We define a TTN as a balanced binary rooted tree of order-three tensors. 
Nodes in the tree are labeled by $t \in \mathcal P(\{1, \ldots, d\})$, 
and each non-leaf node has two childen $t_\mathrm{L}$ and $t_\mathrm{R}$, forming an ordered partition of label $t$.
The root node or base node has label $\tb = D = \{1, \ldots, d\}$ and 
leaf or external nodes have labels $t \in \{\{1\}, \ldots, \{d\}\}$ associated to \emph{physical dimensions} $n_1, \ldots, n_d \in \N$.
The set of all external nodes is denoted by $E$, whereas internal nodes belong to set $J$. 
Internal nodes adjacent to a node in $E$ are furthermore collected in set $L \subset J$.
It is sometimes convenient to consider the root node $\tb$ not as internal but as external. 
Then we will write $J^- = J  \setminus \tb$ and $E^+ = E \cup \tb$. 
The set of all nodes is denoted by $T \subset \mathcal P(\{1, \ldots, d\})$.

\emph{Transfer tensors \& $t$-frames.} To each internal node, we assign a \emph{transfer tensor} $\mathbf{B}_t \in \R^{k_{\tl} \times k_{\tr} \times k_t}$, 
where $\mathbf k = (k_t)_{t \in J}$ is the tuple of \emph{bond dimensions}. The set of all transfer tensors forms a TTN-\emph{parameter}
\begin{equation*}
    x = (\mathbf{B}_t) = (\mathbf{B}_t)_{t \in J} \in \mathcal E,
\end{equation*}
where $\mathcal E = \varprod_{t \in J} \R^{k_{\tl} \times k_{\tr} \times k_t}$ is a linear space.
Additionally,
any node in the network is associated with a \emph{$t$-frame} $U_t \in \R^{n_t \times k_t}$, 
that fulfills the following important recursion:
\begin{equation}\label{eq:ttn_recursion}
    U_t = (U_{\tl} \otimesk U_{\tr}) B_t \text{.}
\end{equation}
Here $n_t$ evaluates as the 
product of all physical dimensions contained within the respective subtree at node $t$.

\emph{Differences to HT.} An important departure from the classical definition of the HT format, 
is that at external nodes $t \in E$,
we enforce $U_t = I \in \R^{n_t \times n_t}$. 
Consequently, the dimensions $(k_t)_{t \in E}$ 
will coincide with the respective physical dimension $(n_t)_{t \in E}$. 
Another change involves the root node $\tb$, where we allow for the \emph{label dimension} $K = k_{\tb}$ to 
take values larger than one.

\emph{Representation.} TTNs are commonly used as efficient representation of high-dimensional tensors. 
Contracting all transfer tensors in the network according to the tree structure constitutes a smooth map
\begin{equation} \label{eq:ttn_tensor}
    \phi: \mathcal E \to \R^{n_1 \times \ldots \times n_d \times K}, \qquad x \mapsto \mathbf{X},
\end{equation}
and the matricization of $\mathbf{X}$ equates the $t$-frame of the base node via 
$
        \mathbf{X}^{(1, \ldots, d)} = U_{\tb} \text{.}
$

\begin{figure}[h]
\centering

\includegraphics[width=0.9\textwidth]{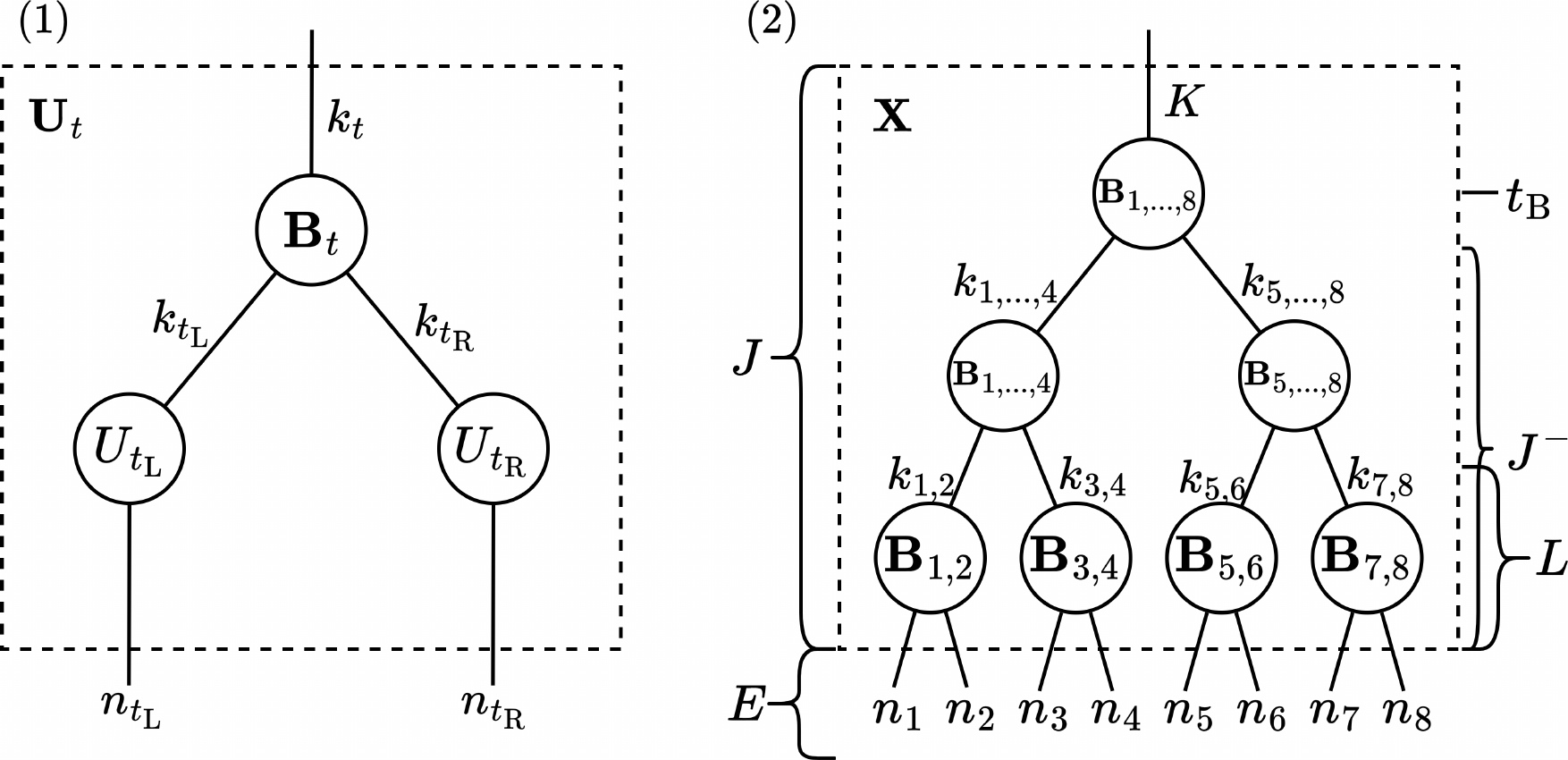} 

\caption{\Cref{eq:ttn_recursion,eq:ttn_tensor} in tensor network diagram notation}
\end{figure}

\subsection{Orthogonal TTNs} \label{sec:orthogonal_ttns}

In a \emph{full-rank} TTN, all transfer tensors $\mathbf{B}_t$ 
are of full rank with respect to matricizations towards inner links of the network \cite[Thm. 1 (ii)]{Uschmajew:2013}. 
The space of full-rank parameters is denoted by $\mathcal{E}^*$, 
and it is an open and dense subset of linear space $\mathcal E$ \cite{Uschmajew:2013}. 
For $\mathcal{E}^*$ to be non-empty, 
the prescribed bond dimensions must satisfy the structural admissibility condition $k_t \leq k_{\tl}k_{\tr}$ for all nodes $t \in J^-$,
such that there exist transfer tensors of full rank. 
We assume throughout that this condition holds.
For a full-rank parameter $x \in \mathcal{E}^*$, the HT-ranks of $\phi(x)$ will coincide 
with the bond dimensions $(k_t)_{t \in J^-}$.

We define the set of \emph{orthogonal} TTN-parameters as
$$\mathcal T = \left\{x = (\mathbf{B}_t) \in \mathcal{E}^* : B_t^\T B_t = I \text{ for all } t \in J^-\right\}.$$
There is no constraint on the root node.
In practice, orthogonal TTN-parameters are preferred over arbitrary parameters for numerical reasons. 
Also in this work, we restrict our analysis to orthogonal TTNs with fixed bond dimensions.
This restriction does not limit 
the expressive power of the TTNs: any TTN parameter $\bar{x}$ can be transformed into an orthogonal parameter $x$ using
HOSVD and reorthogonalization algorithms \cite[Alg. 1 \& 6]{Kressner:2012}, such that $\phi(\bar x) = \phi(x)$. 

The set $\mathcal T$ is an embedded manifold of $\mathcal E$
and its tangent space reads
\begin{equation*} 
\begin{split}
    \TT & = \left(\varprod_{t \in  J^-} T_{B_t}\St(k_{\tl}k_{\tr}, k_t) \right)\times \R^{k_{\tl}k_{\tr} \times K}
    \\&= \left\{ \xi_x = (\delta \mathbf{B}_t):
        B_t^\T \delta B_t  + \delta B_t^\T B_t = 0 \text{ for all } t \in {J^-}
    \right\},
\end{split}
\end{equation*}
where $\delta \mathbf B_t \in \R^{k_{\tl} \times k_{\tr} \times k_t}$ represents a variation of tensor $\mathbf B_t$.
These properties follow from the identification $x = (\mathbf{B}_t) \cong (B_t)$, and noting that
\begin{equation*} 
        \mathcal{T} \cong \left(\varprod_{t \in  J^-} \text{St}(k_{\tl}k_{\tr}, k_t) \times \R^{k_{\tl}k_{\tr} \times K}\right) \cap \mathcal{E}^*
    \end{equation*}
with $\text{St}(m, k) \subset \R^{m\times k}$ ($k \leq m$) the Stiefel manifold. It is well known that Stiefel manifolds 
    are embedded submanifolds \cite[Sec. 3.3.1]{Absil:2008}. Therefore, $\mathcal{T}$ is an embedded submanifold of $\mathcal{E}$ as a 
    Cartesian product of embedded submanifolds, 
    and its tangent space is a linear subspace of $\mathcal E$,
    composed as a Cartesian product of Stiefel 
    tangent spaces \cite[Prop. 5.37]{Lee:2000}. Intersecting with $\mathcal{E}^*$ does not alter those properties, 
    because it is an open subset and therefore intersects the Cartesian product transversely as a special case of \cite[Thm. 6.30]{Lee:2000}. 

    On $\mathcal{T}$, we will work with the Euclidean metric inherited from its linear embedding, which for some $x \in \mathcal{T}$ and $\xi_x, \eta_x \in \TT$ we can write as
    \begin{equation*}
        g_x(\xi_x, \eta_x) = \dotpro{\xi_x}{\eta_x} = \sum_{t \in J} 
        \dotpro{\delta \mathbf{B}_t}{\delta \mathbf{C}_t} .
    \end{equation*}
    With this, $\mathcal{T}$ trivially is a Riemannian submanifold of $\mathcal E$.
    The image of $\mathcal T$ (and also $\mathcal E^*$) under $\phi$ 
    forms another smooth embedded manifold $\mathcal{H}$, the \emph{fixed-rank HT manifold} \cite[Thm. 4]{Uschmajew:2013}.

\subsection{Quotient structure}
It is well-established that full-rank parameters possess a quotient structure with respect to map $\phi$ \cite[Sec. 4.1]{Uschmajew:2013}.
In this subsection, we specialize those results to orthogonal TTN-parameters. 

Define the Lie group $\mathcal{G} = \varprod_{t \in J^-} \Ortho(k_t)$
and the smooth, proper and free Lie group action 
\begin{equation*} 
\begin{split}
    &\theta: \mathcal{T} \times \mathcal{G} \to \mathcal{T}, \\
    & \underbrace{((\mathbf{B}_t), (A_t)_{t \in J^-})}_{= (x, \mathcal{A})} \mapsto  
    \underbrace{((A_{\tl}^\T \otimes A_{\tr}^\T \otimes A_t^\T)\mathbf{B}_t)}_{= \theta(x, \mathcal{A})},
    \quad \text{with } A_t = I \text{ for } t \in E^+.
\end{split}
\end{equation*}
Frequently, we will fix $\mathcal{A} \in \mathcal{G}$, yielding the $g$-isometric diffeomorphism $\theta_\mathcal{A}: \mathcal{T}\to \mathcal{T}$.
Importantly, $\phi$ is invariant under action $\theta$, i.e. $\phi(x) = \phi(\theta(x, \mathcal{A}))$ for all $\mathcal A \in \mathcal G$. By defining the equivalence relation $x \sim \theta(x, \mathcal{A})$, it induces
a quotient manifold 
$
    \mathcal{T}/\mathcal{G} = \{[x]: x \in \mathcal{T}\}, 
$
where $[x] = \{y \in \mathcal{T}: x \sim y\}$. We denote the respective quotient map by
\begin{equation*}
    \pi: \mathcal{T} \to \mathcal{T / G}, \qquad x\mapsto [x] ,
\end{equation*}
and $\mathcal G_x = \pi^{-1}([x])$, the fiber over $[x]$, is an embedded submanifold of $\mathcal{T}$. Pushing $\phi$ down to the quotient yields the map $\hat{\phi}: \mathcal{T/G} \to \R^{n_1 \times \ldots \times n_d \times K}$, which 
is injective as opposed to $\phi$.

Note that whenever we want to provide statements that apply more generally than just to $\mathcal{T}$, the letter $\mathcal{M}$ will always denote some smooth manifold acted upon by a
smooth, free and proper Lie group action $\theta$, and $\mathcal{M/G}$ will denote the 
corresponding quotient manifold for some Lie group $\mathcal{G}$.

\subsection{Quotient optimization problem}
In this work, we want to solve the following problem for some smooth function $h:\R^{n_1 \times \ldots \times n_d \times K} \to \R$ 
using iterative optimization algorithms:
\begin{equation} \label{eq:optimization_problem}
   \min_{\mathbf X \in \mathcal H} h(\mathbf X) = \min_{x \in \mathcal{T}} h (\phi(x)) = \min_{[x] \in \mathcal{T/G}} h(\hat{\phi}([x])).
\end{equation}
We assume that the minimizer of $h$ in $\R^{n_1 \times \ldots \times n_d \times K}$ either admits a representation in TTN format with low ranks or is approximable by such.
If it does, 
limiting our search space to $\mathcal{T}$ is not only justifiable but advisable. 
The bond dimensions $(k_t)_{t \in J^-}$ remain fixed throughout our considerations and
can be seen as an optimization hyperparameter, 
or might be chosen dynamically.

Ideally, we would like to work with the right-hand side of \eqref{eq:optimization_problem} due to the injectiveness of $\hat{\phi}$, 
but elements and tangent vectors of the quotient are abstract objects that are hard to work with numerically. 
This motivates fixing a \emph{horizontal distribution} on the parameter space and leveraging it for optimization.
While optimization algorithms themselves will be formulated on the quotient, the actual implementation happens through 
invariant lifts to the total space $\mathcal{T}$. The main goal of this work will be to assess how different horizontal distributions
affect both theoretical and numerical aspects of quotient optimization on tree tensor networks.

\section{Horizontal Spaces}
\label{sec:manifold} 
We begin this section by establishing a framework that is versatile enough
to allow for arbitrary horizontal distributions.

\subsection{Vertical and horizontal spaces}

\begin{definition}
\label{def:vertical_and_horizontal}
Let $\mathcal{M}$ be a smooth manifold and let $\mathcal{M/G}$ be its quotient. The \emph{vertical space} $\VM$ at $x \in \mathcal{M}$ is the subspace
\begin{equation*} 
    \VM = \ker \dd \pi(x) = T_x\mathcal{G}_x.
\end{equation*}
A \emph{horizontal distribution} is a smooth choice of complementary subspaces to the vertical spaces $\VM$, 
such that any such \emph{horizontal space} $\HM$ 
is invariant with respect to the Lie group action $\theta$, i.e.
\begin{equation} \label{eq:horizontal_compatibillity}
        \dd \theta_{\mathcal{A}}(x)[\HM] = \HM[\theta(x, \mathcal{A})] \text{ for all } x \in \mathcal{M} \text{ and } \mathcal{A} \in \mathcal{G}.
\end{equation}

\end{definition}
\begin{figure}[h]
    \centering
    \includegraphics[width=0.7\textwidth]{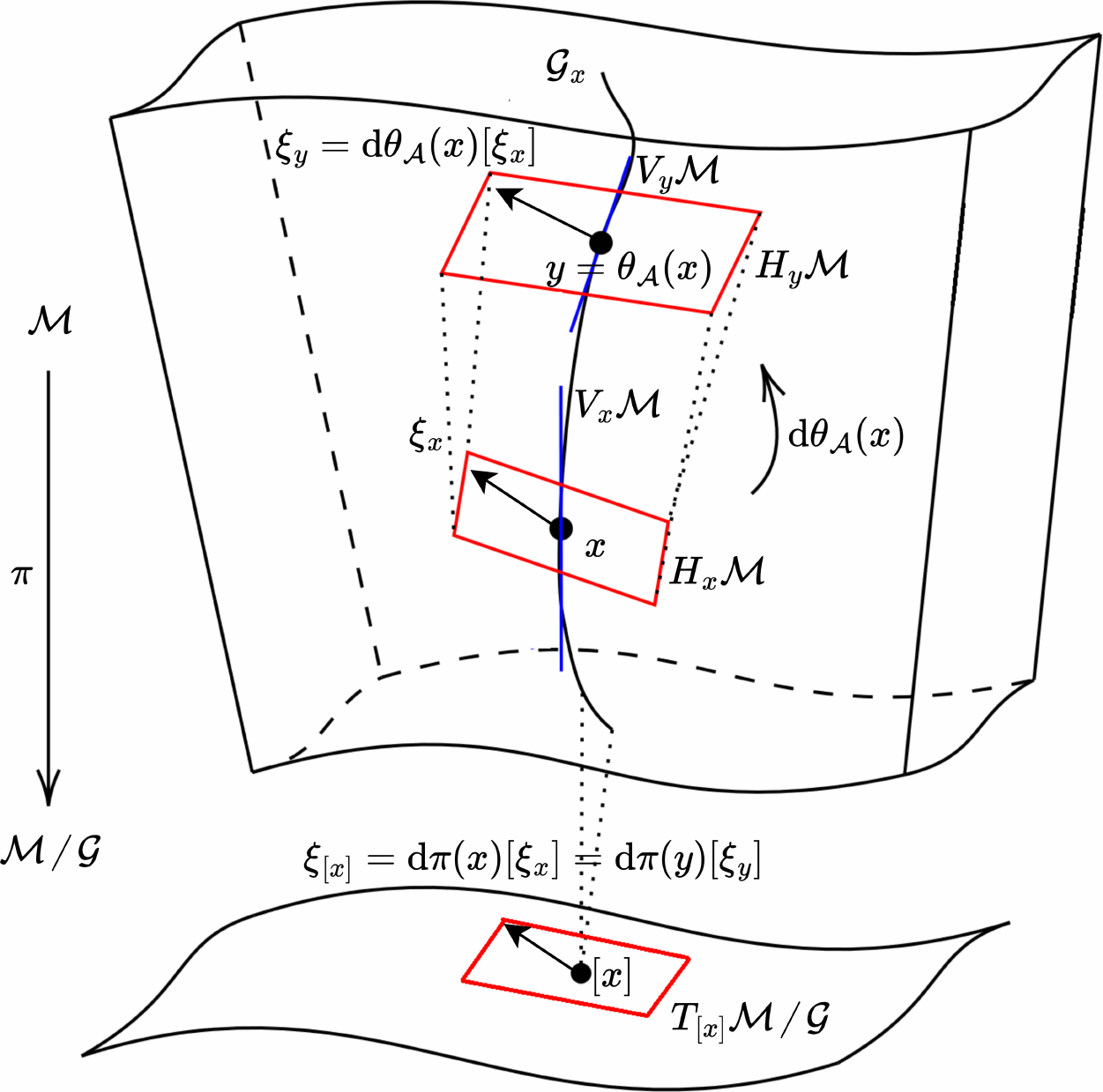}
    \caption{The total space $\mathcal{M}$ and the quotient space $\mathcal{M/G}$. 
    Both $x$ and $y \in \mathcal{G}_x$ map to the same point $[x]$ under the quotient map $\pi$. Vertical spaces run tangent to the orbit $\mathcal{G}_x$. 
    Horizontal spaces along $\mathcal{G}_x$ are compatible with $\dd \theta_\mathcal{A}$. $\xi_x$ and $\xi_y$ are the 
    unique horizontal lifts of $\xi_{[x]}$ to $x$ and $y$, respectively. Adapted from \cite{Wikimedia01}. }
    \label{fig:tangent_space}
\end{figure}

The restriction of $\dd \pi(x)$ to $\HM$ yields a linear isomorphism that
allows the definition of the \emph{horizontal lift} for $\xi_{[x]} \in \TMG$:
\begin{equation*} 
    \lift_x(\xi_{[x]}) = \left(\dd \pi(x) |_{\HM}\right)^{-1}[\xi_{[x]}].
\end{equation*}
In contrast to many modern treatments such as \cite[Def. 9.24]{Boumal:2023}, \cite[Prop. 2.25]{Lee:2019}, \cite[Def. 18.2]{Gallier:2020},
here we do not require $\HM$ to be the orthogonal complement of $\VM$ in $\TM$. Instead we opt for the more general definition of \cite[Def. 10.1]{Nakahara:2003},
which only demands compatibility with $\dd \theta$, as a substitute for orthogonality. The \emph{orthogonal horizontal space} is still included in this definition,
as seen in the next proposition.

\begin{proposition}
    Let ($\mathcal{M}, g$) be a Riemannian manifold, let $\theta$ be a Riemannian isometry on $\mathcal{M}$ and 
    let $\HO$ be the orthogonal complement of $\VM$ in $\TM$. Then $\HO$ is a horizontal space of $\mathcal{M}$ with respect to Lie group $\mathcal{G}$.
\end{proposition}
\begin{proof}
    Let $x \in \mathcal{M}, y = \theta_{\mathcal{A}}(x)$ and $\xi_x \in \HO$. Note that $\theta_\mathcal{A}(\mathcal{G}_x) = \mathcal{G}_y$, 
    and as $\theta_\mathcal{A}$ is diffeomorphic, it holds that $\dd\theta_\mathcal{A}(x)[T_x\mathcal{G}_x] = T_y\mathcal{G}_y$. 
    Using that $\theta$ is an isometry, we have
    \begin{equation*}
        g_x(\xi_x, T_x\mathcal{G}_x) = 0 \Rightarrow g_y(\dd\theta_\mathcal{A}(x)[\xi_x], 
        T_y\mathcal{G}_y) = 0 \Rightarrow \dd\theta_\mathcal{A}(x)[\xi_x] \in \HO[y] .
    \end{equation*}
    As $\theta_\mathcal{A}$ is a diffeomorphism, this already allows to conclude.
\end{proof}

Note that as of \cite[Prop. II.1.2]{Kobayashi:1963}, 
we can equivalently express the invariance \eqref{eq:horizontal_compatibillity} 
of some horizontal space by
\begin{equation} \label{eq:theta_invariance}
        \dd\theta_{\mathcal{A}}(x)[\lift_x(\xi_{[x]})] = \lift_{\theta(x, \mathcal{A})}(\xi_{[x]}) \text{ for all } [x] \in \mathcal{M/G} \text{ and } \xi_{[x]} \in \TMG.
    \end{equation}
Invariance is an important property of horizontal spaces, as it allows sending differential constructs to the quotient, 
by invariantly lifting involved operands to the total space. For example, this can be done
for the Riemannian metric.

\begin{proposition} \label{prop:quotient_metric}
    Let $(\mathcal{M}, g)$ be a Riemannian manifold and let $\theta$ be isometric on $\mathcal{M}$. 
    Any horizontal space $\HM$ induces a Riemannian metric $\hat{g}$ on $\mathcal{M/G}$ via
    \begin{equation} \label{eq:quotient_metric}
        \hat{g}_{[x]}(\xi_{[x]}, \eta_{[x]}) =  g_x(\lift_x(\xi_{[x]}),\lift_x(\eta_{[x]})).
    \end{equation}
\end{proposition}
\begin{proof}
    Let $x \in \mathcal{M}, y = \theta(x, \mathcal{A})$. 
    Employing that $\theta$ is an isometry and \eqref{eq:theta_invariance} delivers
    \begin{equation*} 
        g_x(\lift_x(\xi_{[x]}), \lift_x(\eta_{[x]})) = g_y(\lift_y(\xi_{[x]}), \lift_y(\eta_{[x]})).
    \end{equation*}
    Thus the proposed metric is invariant of the representative in the total space, 
    which qualifies it as well-defined quotient metric \cite[Thm. 9.35]{Boumal:2023}.
\end{proof}

\subsection{Horizontal spaces of TTNs} 
\label{sec:horizontal_spaces}
This section will be dedicated to presenting explicit expressions for vertical and horizontal spaces for the manifolds of orthogonal TTN-parameters.
While all developments of this section and upcoming \cref{sec:projectors} could be extended to full-rank parameters $\mathcal{E}^*$ using
a quotient over $\mathcal{G} = \varprod_{t \in J^-} \mathrm{GL}(k_t)$, this would result in significantly more complex expressions and ultimately to numerically less robust 
algorithms. Thus we limit our exposition to orthogonal parameters $\mathcal T$, for which the vertical space reads
\begin{equation*} 
\VT = 
\left\{\xi_x = (\delta \mathbf{B}_t):
    \begin{aligned}
        \delta \mathbf{B}_t = -G_{\tl} \times_1 \mathbf{B}_t - G_{\tr} \times_2 \mathbf{B}_t + G_{t}^\T \times_3 \mathbf{B}_t \text{ for all } &{t}\in{J}  \\
        G_t \in \Skew(k_t) \text{ for all } t \in J^-  \text{ and } G_{t} = 0 \text{ for all } &t \in E^+
    \end{aligned}
\right\}.
\end{equation*}
It can be found as restriction of \cite[Eq. 25]{Uschmajew:2013}, where it was developed for $\mathcal E^*$.
For the horizontal space, we get two plausible candidates.
The first option arises as a Cartesian product of Stiefel horizontal spaces $H_X\St(m, k)$, 
which is why we call it the \emph{Cartesian horizontal space} of TTNs:
\begin{equation*} 
    \HSt = \left\{\xi_x = (\delta \mathbf{B}_t):
        B_t^\T \delta B_t = 0_{k_t}  \text{ for all } {t}\in{J^-}
     \right\}.
\end{equation*}
It was first proposed by \cite{Holtz:2010}, and is a well-established standard for Riemannian optimization on low-rank tensor train manifolds, because
it provides a computationally feasible parametrization of $T_\mathbf{X}\mathcal H$ (or rather its tensor train analogue).
In \cite{Uschmajew:2013}, it was first used 
as a tool in the quotient geometry of full-rank TTN parameters, and $\HSt$ represents a specialization to 
orthogonal parameters. 
Even though the Cartesian horizontal space is invariant under $\dd \theta$ \cite[Prop. 5]{Uschmajew:2013}, 
it is not orthogonal to $\VT$ with respect to the Euclidean metric $g$, which stands in contrast to 
developments in the literature for Riemannian quotient optimization on matrix manifolds \cite{Absil:2008, Boumal:2023}, 
where the orthogonal horizontal space is considered standard, 
and contemplations on non-orthogonal horizontal spaces remain scarce. This discussion gives rise to our second option.

\begin{theorem}
    The orthogonal complement of $\VT$ in $\TT$ is given by
\begin{equation*} 
    \HT = \left\{ (\delta \mathbf{B}_t) \in \TT:
    \begin{aligned}
        &\delta B_{\tl}^\T B_{\tl} - \mat{B}{1}_{t} (\delta\mat{B}{1}_{t})^\T \in \mathrm{Sym}(k_\tl)\\
        &\delta B_{\tr}^\T B_{\tr} - \mat{B}{2}_{t} (\delta\mat{B}{2}_{t})^\T \in \mathrm{Sym}(k_\tr)
    \end{aligned} \text{ for all $t \in J \setminus L$}
     \right\}.
\end{equation*}
\end{theorem}

\begin{proof}
    We first establish orthogonality. For this, we introduce a new notation for calculating the metric on subtrees of $T$. 
    Let $T^t = \{ s \in T: s \subseteq t\}$ be a subtree of $T$, and let $\xi_x = (\delta \mathbf{B}_t), \eta_x = (\delta \mathbf{C}_t)\in \TT$. Then define
    \begin{equation*}
        \dotpro{\xi_x}{\eta_x}^t = \sum_{s \in T^t} \dotpro{\delta \mathbf{B}_s}{\delta \mathbf{C}_s} .
    \end{equation*}
    It follows that
    \begin{equation*}
    \begin{split}
        & \dotpro{\xi_x}{\eta_x}^t = \dotpro{\xi_x}{\eta_x}^{{\tl}}  + \dotpro{\xi_x}{\eta_x}^{{\tr}} + 
        \dotpro{\delta \mathbf{B}_t}{\delta \mathbf{C}_t} \text{ for all } t \in J \setminus L, 
        \text{ and}\\
        & \dotpro{\xi_x}{\eta_x}^t = \dotpro{\delta  \mathbf{B}_t}{\delta \mathbf{C}_t  } \text{ for all } t \in L .
    \end{split}
    \end{equation*}
    Now let $\xi_x \in \HT$ and $\eta_x \in \VT$. We know we can write $\delta C_t = -(G_{\tl} \otimesk I)B_t - (I \otimesk G_{\tr})B_t + B_t G_t$. 
    Assume that $\dotpro{\xi_x}{\eta_x}^s = \dotpro{\delta B_s}{B_s G_s}$ holds for $s \in \{\tl, \tr\}$. 
    Then utilizing the cyclic property of the trace, we obtain
    \begin{equation*} \label{eq:long_calculation}
    \begin{split}
        \dotpro{\xi_x}{\eta_x}^t =\, & \dotpro{\xi_x}{\eta_x}^{{\tl}}  + \dotpro{\xi_x}{\eta_x}^{{\tr}} + \dotpro{\delta \mathbf{B}_t}{\delta \mathbf{C}_t} \\ 
        =\, &  \dotpro{\delta B_{\tl}}{B_{\tl} G_{\tl}} + \dotpro{\delta B_{\tr}}{ B_{\tr} G_{\tr}} \\ 
        &- \scapro{\delta B_t}{G_{\tl} \otimesk I}{B_t} - \scapro{\delta B_t}{I \otimesk G_{\tr}}{B_t} + \dotpro{\delta B_t}{B_t G_t} \\
        =\, & \Tr{(\delta B_{\tl}^\T B_{\tl} - \mat{B}{1}_{t} (\delta\mat{B}{1}_{t})^\T) G_{\tl}} + 
        \Tr{(\delta B_{\tr}^\T B_{\tr} - \mat{B}{2}_{t} (\delta\mat{B}{2}_{t})^\T) G_{\tr}} \\& + 
        \dotpro{\delta B_t}{B_t G_t} \\
        =\, & \dotpro{\delta B_t}{B_t G_t}.
    \end{split}
    \end{equation*}
    The traces vanished as $\Tr{AB} = 0$ for any symmetric matrix $A$ and skew-symmetric matrix $B$. 
    This justifies the definition of $\HT$ provided in the theorem. Thus, by induction, starting on nodes $t \in L$ that satisfy $\delta C_t = B_t G_t$, 
    the assumption $\dotpro{\xi_x}{\eta_x}^t = \dotpro{\delta B_t}{B_t G_t}$ holds for all $t \in J$. We know however that $G_{\tb} = 0$, 
    so we can infer that $\dotpro{\xi_x}{\eta_x} = 0$. 
    
    The complementarity follows from a dimensionality argument. 
    Note how $\HT$ imposes symmetry conditions on all nodes $t \in J^-$, 
    at each node removing $\tfrac{1}{2} \,k_t(k_t-1)$ degrees of freedom. 
    This is just the number of degrees of freedom that $\VT$ provides: At each node $t \in J^-$, we get to pick a skew-symmetric matrix
    of size $k_t \times k_t$.
    Therefore $\dim \HT + \dim \VT = \dim \TT$, 
    which completes the proof.
\end{proof}

Now that we have established two different horizontal spaces, 
we can lift any vector in $\TTG$ to either $\HSt$ or $\HT$. 
To avoid notational conflicts, we denote the respective lift operations as 
$\lift^\cartesian_x$ and $\lift^\times_x$. Additionally, $\theta$ being an isometry, 
we can leverage \Cref{prop:quotient_metric} to induce two separate Riemannian metrics on the quotient. 
We define
\begin{align*}
    \hat{g}^\cartesian_{[x]}(\xi_{[x]}, \eta_{[x]}) &=  \dotpro{\lift^\cartesian_x(\xi_{[x]})}{\lift^\cartesian_x(\eta_{[x]})}, \mathrm{and}\\
    \hat{g}^\times_{[x]}(\xi_{[x]}, \eta_{[x]}) &=  \dotpro{\lift^\times_x(\xi_{[x]})}{\lift^\times_x(\eta_{[x]})}.
\end{align*}

\subsection{Projectors} \label{sec:projectors}
For optimization on quotient manifolds $\mathcal{M/G}$, 
of which the total space $(\mathcal{M}, g)$ is a Riemannian submanifold of a Euclidean space $\mathcal{E}$, 
the following projectors play an important role:
\begin{itemize}
    \item $\Proj_x: \mathcal{E} \to \TM$, the orthogonal projector onto the tangent space of $\mathcal{M}$;
    \item $\Proj^\times_x: \TM \to \HO$, the orthogonal projector onto $\HO$;
    \item $\Proj^H_x: \TM \to \HM$, denoting an orthogonal projector onto a horizontal space;
    \item $\Proj^{H,V}_x: \TM \to \HM$, denoting an oblique projector along the vertical space onto a horizontal space.
\end{itemize}
Notably, $\Proj^H_x$ and $\Proj^{H,V}_x$ may be constructed for any horizontal space $\HM$, 
but they coincide with $\Proj^\times_x$, if $\HM$ is chosen orthogonal.
This section will cover projectors for TTN-manifold $\mathcal{T}$, starting with those onto $\TT$ and $\HSt$.
Having identified both as Cartesian products of $T_X\St(m, k)$ and $H_X\St(m, k)$ respectively, 
we may perform those projections componentwise. 
They are given by \cite[Sec. 9.16]{Boumal:2023}
\begin{equation*} \label{eq:projectors_stiefel}
\begin{split}
    &\Proj_{X}: \R^{m\times k} \to T_X\St(m, k), \qquad \delta X \mapsto \delta X - \frac{1}{2}\, X (X^\T\delta X + \delta X^\T X),\\
    &\Proj^H_{X}: \R^{m\times k} \to H_X\St(m, k), \qquad \delta X  \mapsto (I - X X^\T)\delta X,
\end{split}
\end{equation*}
and we obtain
\begin{align*}
    \Proj_x: \mathcal{E} \to \TT, \qquad &(\delta \mathbf{B}_t, \delta \mathbf{B}_{\tb}) \mapsto 
    ( \Proj_{\mathbf{B}_t}(\delta \mathbf{B}_t), \delta \mathbf{B}_{\tb}) \text{ and}\\
    \Proj_x^\cartesian: \TT \to \HSt, \qquad &(\delta \mathbf{B}_t, \delta \mathbf{B}_{\tb}) \mapsto 
    ( \Proj^H_{\mathbf{B}_t}(\delta \mathbf{B}_t), \delta \mathbf{B}_{\tb}).
\end{align*}
The projections for tensor-valued components are implemented
by $(1,2)$-matricizing prior to the operation and $(1,2)$-dematricizing the result. 
Here $\Proj_x^\cartesian$ can be seen as a concrete implementation of $\Proj_x^H$ for $\HSt$. 

Implementing a projector for $\HT$ is more challenging, and was posed as an open question in \cite{DaSilva:2015}.
A direct construction of $\Proj_x^\times: \TT \to \HT$ results in a system of $|J|$ coupled Lyapunov equations, 
which need to be solved for each projection. 
Here, we present an alternative approach involving $\Proj_x^{H,V}: \TT \to \HSt$, 
the oblique projector along $\VT$ onto $\HSt$, and the Moore-Penrose pseudoinverse.
\begin{theorem}
    The orthogonal projector onto $\HT$ is given by
\begin{equation} \label{eq:proj_H*}
    \Proj_x^\times = (\Proj^{H,V}_x)^*\circ(\Proj_x^{H,V} \circ (\Proj_x^{H,V})^*)^+\circ\Proj_x^{H,V}.
\end{equation}
\end{theorem}
The result itself is elementary; the interest lies in the construction of the oblique projector,
which comes in the following form:
\begin{equation} \label{eq:projhv}
    \Proj_x^{H,V}:\, (\delta \mathbf{B}_t) \mapsto 
    \Proj_x^\cartesian(
        G_{\tl} \times_1 \mathbf{B}_t + G_{\tr} \times_2 \mathbf{B}_t + \delta \mathbf{B}_t
    )
\end{equation}
with $G_t$ calculated recursively from the leaves up as
\begin{equation*}
\begin{split} 
    G_t = \begin{cases}
        0 &\text{for } t \in E, \\
        B_t^\T\delta B_t + B_t^\T(G_{\tl} \otimesk I)B_t + B_t^\T(I \otimesk G_{\tr})B_t &\text{for } t \in J.
    \end{cases}
\end{split}
\end{equation*}
Its adjoint $(\Proj_x^{H,V})^*: \TT \to \HT$ reads
\begin{equation} \label{eq:oblique_proj_h*}
    (\Proj_x^{H,V})^*: (\delta \mathbf{B}_t) \mapsto (\delta \mathbf{B}_t - G_t^\T \times_3 \mathbf{B}_t).
\end{equation}
This time $G_t$ is recursively calculated from the root down via
\begin{equation*}
\begin{split}
    G_{\tl} &= \mskew(\mat{B}{1}_t(\delta \mat{B}{1}_t)^\T - \delta B_{\tl}^\T B_{\tl}) + \mat{B}{1}_t ( I \otimesk G_t) (\mat{B}{1}_t )^\T, \\
    G_{\tr} &= \mskew(\mat{B}{2}_t(\delta \mat{B}{2}_t)^\T - \delta B_{\tr}^\T B_{\tr}) + \mat{B}{2}_t ( I \otimesk G_t) (\mat{B}{2}_t )^\T,
\end{split}
\end{equation*}
for $t \in J$, starting with $G_{\tb}=0$. 
We omit an explicit derivation of the above statements, because operators 
\eqref{eq:proj_H*}, \eqref{eq:projhv} and \eqref{eq:oblique_proj_h*} 
are easily verified to be projectors onto the specified images along the corresponding kernels by direct calculation. 
For the two latter, it is instructive to plug in general vectors of $\VT$ and
$\HT$ and employ the properties of involved tensors.

The application of \eqref{eq:proj_H*} remains computationally challenging. 
Due to recursive, node-connecting relations of $G_t$,
the transformation matrices of the oblique projectors are not explicitly available, 
and we cannot perform the required inversion directly. 
The projection is instead computed by solving the linear system associated with the pseudoinverse in $\HSt$.

\section{Optimization Tools} \label{sec:opt_tools}
\subsection{Riemannian gradient}
The goal of this section is to find Riemannian gradients of functions on $\mathcal{T/G}$ by 
computing their horizontal lifts to $\mathcal{T}$. This is possible in general, so again consider $\mathcal{M/G}$, 
the quotient of $(\mathcal{M}, g)$, which in turn is
a Riemannian submanifold of a Euclidean space
$\mathcal{E}$. Let $\hat f: \mathcal{M/G} \to \R$, $f = \hat{f} \circ \pi: \mathcal{M} \to \R$ and let
$\bar f: \mathcal{E} \to \R$ be a smooth extension of $f$ onto the surrounding Euclidean space.
Standard results \cite[Eq. 3.47 \& Eq. 3.39]{Absil:2008} deliver that
\begin{align} 
    \grad f(x) &= \Proj_x(\nabla \Bar{f}(x)), \nonumber \\
    \lift^\times_x(\grad^\times\hat{f}([x])) &= \grad f(x), \label{eq:gradient_lift*}
\end{align}
where $\grad^\times\hat{f}$ is the quotient gradient associated with metric $\hat g^\times$ induced by the orthogonal horizontal space in the sense of \Cref{eq:quotient_metric}.
While the generalization to arbitrary quotient metrics $\hat{g}$ is not explicitly detailed in the literature, 
it follows by a straightforward modification of the proof in \cite[Prop. 9.39]{Boumal:2023}, resulting in
\begin{equation} \label{eq:gradient_lift}
    \lift_x(\grad \Hat{f}([x])) = \Proj^H_x(\grad f(x)).
\end{equation}
Notably, \eqref{eq:gradient_lift} does not recover the stronger result \eqref{eq:gradient_lift*}, 
as $\HM$ lacks the orthogonality to $\VM$ in general.

Recall that we want to solve \eqref{eq:optimization_problem} for some 
smooth function $h: \R^{n_1 \times \ldots \times n_d \times K} \to \R$. 
Define the functions $\Hat{f} = h \circ \Hat{\phi}: \mathcal{T/G} \to \R$ and 
$f = h \circ \phi|_\mathcal{T}: \mathcal{T} \to \R$. Importantly, it holds that $f = \hat{f} \circ \pi$, 
and $f$ admits a natural smooth extension $\bar f = h \circ \phi$ onto $\mathcal{E}$, so all of the above statements apply to $\mathcal{T}$. 
We will write $\grad^\cartesian \Hat{f}$ for the gradient accompanying $\HSt$, 
and $\grad^\times \Hat{f}$ for the one belonging to $\HT$.

The Euclidean gradient of $\bar f$ is given by $\dd\phi(x)^*[\nabla h (\phi(x))] = \nabla \Bar{f}(x)$. 
An expression for $\dd\phi(x)^*: \R^{n_1 \times \ldots \times n_d \times K} \to T_x\mathcal E$ was derived in \cite[Eq. 13]{DaSilva:2015} for the 
Hierarchical Tucker format and it is applicable to our notion of TTNs.
Let $\delta \mathbf{X} \in \R^{n_1 \times \ldots \times n_d \times K}$ and define intermediate variations $(\mathbf W_t)_{t\in J}, \mathbf W_t \in \R^{n_{\tl} \times n_{\tr} \times k_t}$.
Then the components of $d\phi(x)^*[\delta \mathbf{X}] = (\delta \mathbf{B}_t)$ can be calculated recursively 
for all $t \in J$ via 
\begin{align} \label{eq:gradient_recursion}
\begin{split}
    &W_{\tl} = \dotpro{U_{\tr}^\T \times_2 \mathbf W_t}{\mathbf{B}_t}_{(2, 3)}\\ &
    W_{\tr} = \dotpro{U_{\tl}^\T \times_1\mathbf W_t}{\mathbf{B}_t}_{(1, 3)} \\ &
    \delta \mathbf{B}_t = (U_{\tl}^\T \otimes U_{\tr}^\T \otimes I_{k_t}) \mathbf W_t
\end{split}
\end{align} 
starting from $W_{\tb} = (\delta \mathbf{X})^{(1,\ldots, d)}$.
The complete chain of transformations required for TTN quotient gradients is depicted in the schematic below,
which we will apply to a machine learning problem in \Cref{sec:ml_ttn}.
\begin{equation*}
\begin{tikzcd}[column sep=large]
\underset{\nabla h}{\R^{n_1 \times \ldots \times n_d \times K }} & \underset{\nabla \bar{f}}{T_x\mathcal{E}} & \underset{\grad f}{\TT} & H_x\mathcal{T} & \underset{\grad \hat f}{\TTG} 
\arrow[from=1-1, to=1-2, swap, "\dd \phi(x)^*"]
\arrow[from=1-2, to=1-3, swap, "\Proj_x"]
\arrow[from=1-3, to=1-4, swap, "\Proj_x^H"]
\arrow[from=1-5, to=1-4, "\lift_x"]
\end{tikzcd}
\end{equation*}

\subsection{Riemannian connection} \label{sec:connection}
To analyze derivatives of vector fields, we will employ the notion of connections described in \cite[Sec. 4]{Lee:2019},
so a connection is a smooth map $\nabla: \mathfrak{X}(\mathcal{M}) \times \mathfrak{X}(\mathcal{M}) \to \mathfrak{X}(\mathcal{M})$, 
that is $\mathrm C^\infty(\mathcal M)$-linear in its first argument, $\R$-linear in its second argument, and that fulfills the Leibniz rule of connections.
$\mathfrak{X}(\mathcal{M})$ denotes the set of all smooth vector fields on a manifold $\mathcal{M}$. 
For $X, Y \in \mathfrak{X}(\mathcal{M})$, $\nabla_X Y$ is called the covariant derivative
of $Y$ in the direction $X$.
Note that standard operations extend to vector fields. 
Firstly, any of the projectors covered in \Cref{sec:projectors} map pointwise between smooth vector fields, 
e.g. $\Proj: \mathfrak{X}(\mathcal{E}) \to \mathfrak{X}(\mathcal{M})$, such that $\Proj(U)(x) = \Proj_x(U(x))$ for all 
$x \in \mathcal{M}, U \in \mathfrak{X}(\mathcal{E})$. Secondly, any horizontal lift is a map 
$\lift:\mathfrak{X}(\mathcal{M/G}) \to \mathfrak{X}(\mathcal{M})$, such that $\lift(U)(x) = \lift_x(U(\pi(x)))$ for all 
$x \in \mathcal{M}, U \in \mathfrak{X}(\mathcal{M/G})$ \cite[Thm. 3]{Uschmajew:2013}. 

Like previously for the Riemannian gradients of $\mathcal{M/G}$, different metrics induce different connections. 
For the metric induced by the orthogonal horizontal space, in the literature \cite{Absil:2008,Boumal:2023} 
there exist constructions of 
a Riemannian connection (i.e. a connection which is torsion-free and compatible with the metric \cite[Sec. 5]{Lee:2019}) on the quotient.
The following theorem, on the other hand, provides the construction of a connection with torsion for any horizontal metric.
\begin{theorem} \label{thm:quotient_connection}
    Let $(\mathcal{M}, g, \nabla)$ be a Riemannian manifold with its Riemannian connection. Then
    \begin{equation} \label{eq:quotient_connection}
        \lift(\hat{\nabla}_U V) = \Proj^H(\nabla_{\lift(U)}\lift(V))
    \end{equation}
    defines a connection $\hat{\nabla}$ on $\mathcal{M/G}$ 
    for any horizontal space $\HM$ and all vector fields $U, V \in \mathfrak{X}(\mathcal{M/G})$.  If $\HM$ induces metric $\hat{g}$ on the quotient, 
    $\hat{\nabla}$ is compatible with that metric, but it is in general not torsion-free. Its torsion tensor reads
    \begin{equation} \label{eq:torsion_tensor}
        \lift(T(U, V)) = \Proj^H([\lift(U), \lift(V)]) - \Proj^{H, V}([\lift(U), \lift(V)]).
    \end{equation}
\end{theorem}
\begin{proof}
    By pushing \eqref{eq:quotient_connection} to the quotient, i.e. applying $\dd \pi$, 
    it becomes apparent that $\hat\nabla$ produces vector fields of $\mathcal{M/G}$.
    Smoothness, linearity and Leibniz rule are inherited by $\nabla$ of $\mathcal{M}$, 
    due to the linearity of $\dd \pi$ and $\Proj^H$, so it is clear that $\hat\nabla$ is a connection.

    For the compatibility with metric $\hat g$ consider $U, V, W \in \mathfrak{X}(\mathcal{M/G})$ 
    and their horizontal lifts $\bar U = \lift(U), \bar V = \lift(V), \bar W = \lift(W)$. Then 
\begin{equation*}
\begin{split}
    {U}\hat{g}(V,W) \circ \pi &= \Bar{U} g(\Bar{V}, \Bar{W}) \\&= g({\nabla_{\Bar{U}} \Bar{V}},{\Bar{W}}) + g({\Bar{V}},\nabla_{\Bar{U}} \Bar{W})\\
    &= g({\Proj^H(\nabla_{\Bar{U}} \Bar{V})},{\Bar{W}}) + g({\Bar{V}},{\Proj^H(\nabla_{\Bar{U}} \Bar{W})}) \\
    &= g({\lift(\Hat{\nabla}_{{U} }{V})},{\lift({W})}) + g({\lift({V})},{\lift(\Hat{\nabla}_{{U}}{W})}) \\
    &=(\hat{g}(\Hat{\nabla}_{{U}}{V}, {W}) + \hat{g}({V},\Hat{\nabla}_{{U}}{W})) \circ \pi.
\end{split}
\end{equation*}
In the first equality, we employed \cite[Eq. 9.46]{Boumal:2023}, in the second that $\nabla$ is compatible with $g$, 
and in the third, we applied $\Proj^H$ to $\Bar{W}$ and $\Bar{V}$ respectively, as they are horizontal anyway, 
and then made use of the orthogonality of $\Proj^H$.

For the torsion, we first need to make some preparations. 
Consider a function $\hat{f}:\mathcal{M/G} \to  \R$ and its lift $f = \hat f \circ \pi$. Then \cite[Eq. 9.26]{Boumal:2023} tells us 
that $(U\hat f) \circ \pi = \bar Uf$. Furthermore, by \cite[Prop. II.1.3]{Kobayashi:1963} it holds that
\begin{equation*}
    \lift([U, V]) = \Proj^{H, V}([\bar U, \bar V]).
\end{equation*}
With this in place we can argue
\begin{equation*}
\begin{split}
    \lift(T(U, V))f &= (T(U, V)\hat f) \circ \pi = ((\hat\nabla_U V - \hat\nabla_V U - [U, V])\hat f) \circ \pi \\
    &= \lift(\hat\nabla_U V)f - \lift(\hat\nabla_V U)f - \lift([U, V])f \\
    &= \Proj^H(\nabla_{\Bar{U}} \Bar{V})f - \Proj^H(\nabla_{\Bar{V}} \Bar{U}) - \Proj^{H, V}([\bar U, \bar V])f \\
    &= (\Proj^{H}([\bar U, \bar V]) - \Proj^{H, V}([\bar U, \bar V]))f.
\end{split}
\end{equation*}
\end{proof}
For the orthogonal horizontal space $\HT$, the above theorem allows us to construct 
the connection $\hat\nabla^\times$, which is even the Riemannian one with respect to metric $\hat g^\times$, 
because the projectors of \eqref{eq:torsion_tensor} coincide and the torsion vanishes identically. 
This fact can also be found in \cite[Prop. 5.5.3]{Absil:2008} or \cite[Thm. 9.43]{Boumal:2023}, 
from which \Cref{thm:quotient_connection} drew inspiration. 
For Cartesian horizontal space $\HSt$ and metric $\hat g^\cartesian$, an expression for the Riemannian connection eludes us so far;
we instead utilize the connection $\hat\nabla^\cartesian$ provided by \Cref{thm:quotient_connection}.

\begin{remark}
    The connection of \Cref{thm:quotient_connection} can also be motivated with the following argument: the map $\mathrm{V}_{[x] \leftarrow [y]}: \TMG[y] \to \TMG[x]$ given by
    \begin{equation*}
        \mathrm{V}_{[x] \leftarrow [y]} = \dd \pi(x) \circ \Proj_x^H \circ \lift_y
    \end{equation*}
    is a numerically viable vector transport on $\mathcal{M/G}$. Then according to \cite[Lemma 2.2]{Weigl:2024},
    a numerically viable connection is given by differentiating this transporter. 
    \Cref{eq:quotient_connection} may be derived in exactly this way.

\end{remark}

\subsection{Covariant Hessian}
\begin{definition} Let $(\mathcal{M}, g, \nabla)$ be a smooth Riemannian manifold and let connection $\nabla$ be compatible with $g$. 
    The \emph{covariant Hessian} of function $f: \mathcal{M} \to \R$ at $x \in \mathcal{M}$ is the linear map $\Hess f_x: \TM \to \TM$ defined as:
\begin{equation*} 
    \Hess f_x ( \xi_x) = \nabla_{\xi_x}(\grad f).
\end{equation*}
\end{definition}
Formulating covariant Hessians of a quotient $\mathcal{M/G}$ 
in terms of horizontal lifts to the total space $\mathcal{M}$ will be the topic of this section. 
Consider a function $\hat{f}:\mathcal{M/G} \to \R$ and let $\Hess \hat{f}$ be the covariant Hessian 
\begin{equation*}
    \Hess \hat{f} (U) = \hat{\nabla}_U(\grad \hat{f})
\end{equation*}
for any connection $\hat\nabla$ delivered by \Cref{thm:quotient_connection}. 
The horizontal lift of this Hessian constitutes a linear map 
\begin{equation*}
    \mathrm{H}_x: \TM \to \TM, \qquad \xi_x \mapsto \nabla_{\xi_x}\lift(\grad \hat{f}),
\end{equation*}
which relates to the quotient via
\begin{equation} \label{eq:hessian_lift}
    \lift_x(\Hess \hat f_{[x]} (\xi_{[x]})) = \Proj^H_x(\mathrm{H}_x(\lift_x(\xi_{[x]}))).
\end{equation}
A generalization of \cite[Ex. 9.47]{Boumal:2023} characterizes the eigenvalues of lifted Hessians of this form.
Specifically, $ \Proj^H_x \circ \mathrm{H}_x \circ \Proj_x^H $ recovers the spectrum of $\Hess \hat{f}_{[x]}$ on $\HM$.

When restricting ourselves to $\HT$ and its Riemannian connection $\hat\nabla^\times$,
the respective covariant Hessian defined via
 \begin{equation*}
     \Hess^\times \hat{f} (U) = \hat\nabla^\times_U(\grad^\times \hat f)
 \end{equation*}
is called the Riemannian Hessian of $\hat{f}$ with respect to $\hat{g}^\times$. 
Furthermore, since in our scenario $f = \hat f \circ \pi$, we can calculate
\begin{align*}
    \lift^\times (\Hess^\times \hat{f} (U)) &\overset{\eqref{eq:quotient_connection}}{=} \Proj^\times (\nabla_{\lift^\times (U)}\lift^\times(\grad^\times \hat{f})) \\
    &\overset{\eqref{eq:gradient_lift*}}{=} \Proj^\times (\nabla_{\lift^\times (U)} \grad f) = \Proj^\times (\Hess f (\lift^\times (U))),
\end{align*} 
to find that the lifted Hessian actually closely relates to the Riemannian Hessian of $\mathcal{T}$ \cite[Prop. 9.45]{Boumal:2023}. This does not hold for $\HSt$: neither the covariant Hessian
 \begin{equation*}
     \Hess^\cartesian \hat{f} (U) = \hat\nabla^\cartesian_U(\grad^\cartesian \hat f)
 \end{equation*}
nor its lift $\mathrm{H}^\cartesian_x$ are Riemannian Hessians.

We will close this section with the following observation:
As seen in the proof of \cite[Prop. 5.15]{Boumal:2023} or equivalently in \cite[Ex. 4-6c]{Lee:2019}, for a connection $\nabla$, 
that is compatible with metric $g$, its covariant Hessian is symmetric if and only if $\nabla$ is torsion-free. 
This means Riemannian Hessians are symmetric maps, whereas $\Hess \hat{f}_{[x]}$ and its lift $\mathrm{H}_x$ are in general non-symmetric.

\subsection{Retractions}
A \emph{retraction} $\Ret: T\mathcal{M} \to \mathcal{M}: (x, \xi_x) \mapsto \Ret_x(\xi_x)$
is a smooth generator of curves $\gamma: [a, b] \to \mathcal M$ with $\gamma(s) = \Ret_x(s\xi_x)$, such that $\gamma(0) = x$ and $\gamma'(0) = \xi_x$ \cite[Def. 3.47]{Boumal:2023}.
It is used in iterative optimization on manifolds to update according to directions, that reside
in the tangent bundle $T\mathcal{M} = \bigcup_{x \in \mathcal M} T_x\mathcal M$ without leaving the manifold itself.
In this section, we will cover several retractions on $\mathcal{T}$ and send them to $\mathcal{T/G}$ using \cite[Prop 4.1.3]{Absil:2008}. 
This proposition importantly states that any retraction $\Ret$ on $\mathcal{T}$ provides a retraction $\hat \Ret$ on $\mathcal{T/G}$ by
\begin{equation} \label{eq:quotient_retraction}
    \hat{\Ret}_{[x]}(\xi_{[x]}) = \pi(\Ret_x(\lift_x(\xi_{[x]}))),
\end{equation}
if it satisfies the following equivariance condition for all $x \sim y \in \mathcal{T}$ and $\xi_{[x]} \in \TTG$:
\begin{equation} \label{eq:retraction_horizontal_invariance}
    \Ret_x(\lift_x(\xi_{[x]})) \sim \Ret_y(\lift_y(\xi_{[x]})).
\end{equation} 
Consequently, retractions on the quotient can be performed implicitly by retracting on the total space.
Note that even though this result was originally formulated for lifts to the orthogonal horizontal space,
it applies to arbitrary horizontal spaces without modification.
\begin{algorithm}[H]
\caption{Cartesian retraction $\Ret^{\cartesian}: (x, \xi_x) \mapsto \Ret^{\cartesian}_x(\xi_x)$}
\label{alg:cartesian_retraction}
\begin{algorithmic}
    \REQUIRE $x = (\mathbf{B}_t) \in \mathcal{T}$, $\xi_x = (\delta \mathbf{B}_t) \in \TT$, retraction $\Ret^{\text{St}}$ on the Stiefel manifold.
    \FOR{$t \in J^-$}
    \STATE $\mat{C}{1,2}_t \gets \Ret^{\text{St}}_{B_t}(\delta B_t)$
    \ENDFOR
    \STATE $\mathbf{C}_{\tb} \gets \mathbf{B}_{\tb} + \delta \mathbf{B}_{\tb}$    
    \RETURN $(\mathbf{C}_t)\in \mathcal{M}$
\end{algorithmic}
\end{algorithm}
Due to $\mathcal{T}$ being a Cartesian product of Stiefel manifolds, retracting componentwisely delivers a retraction. Cartesian retractions are computationally very attractive, because they act on each tensor separately, and can easily be parallelized, as opposed 
to the recursive retractions presented in \cite[Alg. 3]{DaSilva:2015} or \cite[Prop. 4.11]{Steinlechner:2016}.
\begin{theorem}
\label{thm:cartesian_retraction} Let $\Ret^{\text{St}}$ be a retraction on the Stiefel manifold $\St(m, k)$. 
Then the application of \Cref{alg:cartesian_retraction} is a retraction 
$\Ret^{\cartesian}$ on $\mathcal{T}$. If $\Ret^{\text{St}}$ additionally fulfils
\begin{equation} \label{eq:retraction_qlinear}
    \Ret^{\text{St}}_{QXW}(Q \,\delta X\, W) = Q\, \Ret^{\text{St}}_X(\delta X)\, W
\end{equation}
for orthogonal matrices $Q \in \Ortho(n)$ and $ W \in \Ortho(k)$, 
then $\Ret^{\cartesian}$ is equivariant \eqref{eq:retraction_horizontal_invariance} and provides a retraction $\hat{\Ret}^{\cartesian}$ on $\mathcal{T/G}$ through \eqref{eq:quotient_retraction}.
\end{theorem}
\begin{proof}
    For the first part, define the curve $\gamma(s) = \Ret^{\cartesian}_x(s\xi_x)$ for some $x \in \mathcal{T}$ and $\xi_x \in \TT$. This curve is given  by
    \begin{equation*}
        \gamma(s) = (\Ret^{\text{St}}_{B_t}(s\delta B_t), \mathbf{B}_{\tb} + s\delta \mathbf{B}_{\tb}).
    \end{equation*}
    With $\Ret^{\text{St}}$ being a retraction, we immediately find $\gamma(0) = x$ and $\Dot{\gamma}(0) = \xi_x$, 
    so ${\Ret}^{\cartesian}$ is also a retraction. 
    For the second part, we will show that $\Ret^{\cartesian}$ satisfies \eqref{eq:retraction_horizontal_invariance}. 
    Let $y \in \mathcal{G}_x$, with $y = \theta(x, \mathcal{A}) = ((A_{\tl}^\T \otimes A_{\tr}^\T \otimes A_t^\T)\mathbf{B}_t)$. Consider a vector $\xi_{[x]} \in T_{[x]}\mathcal{T/G}$ and its horizontal lifts $\lift_x(\xi_{[x]})$ and $\lift_y(\xi_{[x]})$. Let
    $\lift_x(\xi_{[x]}) = (\delta \mathbf{B}_t)$, then by \eqref{eq:theta_invariance} we have $\lift_y(\xi_{[x]}) = ((A_{\tl}^\T \otimes A_{\tr}^\T \otimes A_t^\T)\delta \mathbf{B}_t)$. Employing the additional condition \eqref{eq:retraction_qlinear} we find
    \begin{equation*}
    \begin{split}
        &\Ret^{\text{St}}_{(A_{\tl} \otimesk A_{\tr})^\T B_t A_t}((A_{\tl} \otimesk A_{\tr})^\T \delta B_t A_t) = (A_{\tl} \otimesk A_{\tr})^\T \Ret^{\text{St}}_{B_t}(\delta B_t) A_t \text{ for all } t \in J^-.
    \end{split}
    \end{equation*}
    Thus, we have shown that $\Ret^{\cartesian}_y(\lift_y(\xi_{[x]})) = \theta(\Ret^{\cartesian}_x(\lift_x(\xi_{[x]})), \mathcal{A})$, 
    so $\Ret^\cartesian$ adheres to \eqref{eq:retraction_horizontal_invariance}
    and $ \Hat{\Ret}^{\cartesian}$ defined like in \eqref{eq:quotient_retraction} is a retraction on $\mathcal{T/G}$.
\end{proof}
Popular retractions for the Stiefel manifold are the \emph{QR decomposition}, the \emph{polar decomposition} \cite[Sec. 7.3]{Boumal:2023}
and the \emph{Cayley transform} \cite{Wen:2010}. Both the polar decomposition and the Cayley transform can be seen to fulfill 
\eqref{eq:retraction_qlinear}, so they are equivariant and deliver retractions on $\mathcal{T/G}$.
The QR decomposition however does not satisfy \eqref{eq:retraction_qlinear}, so we get no well-defined quotient retraction. 
When moving from $[x] \in \mathcal{T/G}$ into a direction $\xi_{[x]}$, 
the resulting path $\hat \gamma: [a, b] \to \mathcal{T/G}$ depends on the chosen representative $y \in \mathcal{G}_x$,
although the initial velocity $\hat \gamma'(0) = \xi_{[x]}$ remains invariant.
From a numerical perspective, this is not detrimental. 
Because we are working with representatives in $\mathcal{T}$ anyway, we can still apply the QR-retraction, 
as verified by the
experiments of \Cref{sec:results}.

\section{Optimization Algorithms} \label{sec:opt_algs}
This section goes along the lines of \cite[Sec. 9.9 \& 9.12]{Boumal:2023}, 
which introduces \emph{Riemannian gradient descent} (RGD),
\emph{Riemannian Newton's method} (RNM) and \emph{Riemannian trust-region} (RTR) for quotient manifolds and 
formulates them in terms of lifts to an orthogonal horizontal space. 
We apply these results to $\mathcal{T}$
and we also consider the case of the non-orthogonal horizontal space $\HSt$.
\subsection{Riemannian gradient descent}
Take an initial point $[x_0] \in \mathcal{T/G}$, some retraction $\hat \Ret$, and suitably chosen step-sizes $\alpha_k$.
Then, for the quotient gradient $\grad \hat{f}$ of some horizontal space, RGD iterates according to
\begin{equation*} 
    [x_{k+1}] = \hat{\Ret}_{[x_k]}(- \alpha_k \grad \hat{f} ([x_k])).
\end{equation*}
Using \eqref{eq:gradient_lift} and \eqref{eq:quotient_retraction}, this iteration rule lifts to the total space as
\begin{align} \label{eq:rgd_H}
    x_{k+1} &= \Ret_{x_k}(-\alpha_k \Proj^\cartesian_{x_k}(\grad f(x_k)))
\end{align}
or alternatively, if $\HT$ is taken as horizontal space, by applying \eqref{eq:gradient_lift*} instead we get
\begin{align} \label{eq:rgd_H*}
    x_{k+1} = \Ret_{x_k}(-\alpha_k \grad f(x_k)).
\end{align}
In practice, one of the two aforementioned update rules is implemented 
in the total space to implicitly optimize on the quotient.
Step sizes $\alpha_k$ can be chosen using Armijo--Goldstein backtracking line search, 
i.e. for some constant $r \in (0, 1)$, $\alpha_k$ is taken such that
\begin{equation*}
    \hat f([x_k]) - \hat f(\hat{\Ret}_{[x_k]}(- \alpha_k \grad \hat{f} ([x_k]))) 
    \geq r\alpha_k \hat{g}_{[x_k]}( \grad\hat{f}([x_k]), \grad \hat{f}([x_k])) .
\end{equation*}
Lifting this condition yields
\begin{equation} \label{eq:armijo_lift}
    f(x) - f(\Ret_{x_k}(-\alpha_k \Proj^\cartesian_{x_k}(\grad f(x_k)))) \geq r \alpha_k ||\Proj^\cartesian_{x_k}(\grad f(x_k))||^2,
\end{equation}
where one might skip the projection, if the orthogonal horizontal space is used.  
When combined with Armijo--Goldstein step-size control, 
local convergence of RGD to a strict second-order critical point can be guaranteed, 
with at least a linear convergence rate \cite[Thm. 4.5.6]{Absil:2008}.
As long as $\hat{f}$ is bounded from below, 
we can usually expect global convergence to a critical point, 
provided that certain regularity conditions are imposed on an equivariant retraction \cite[Corollary  4.13]{Boumal:2023}.

\begin{algorithm} 
\caption{Riemannian gradient descent}
\label{alg:rgd}
\begin{algorithmic}
\REQUIRE Objective function $f: \mathcal{T} \to \R$,
initial TTN-Parameter $x_0 \in \mathcal{T}$, tolerance $\varepsilon$, \\
projector $\mathrm P_x$,
retraction $\mathrm R_x$

\FOR{$k$ = $0, 1, 2, \ldots$}
\STATE $\nabla \Bar{f} \gets$ EuclideanGradient$(f, x_k)$
\STATE $s_k \gets \mathrm{P}_{x_k}(\nabla \Bar{f})$
\IF{$||s_k||^2 < \varepsilon$}
\RETURN $x_k$
\ENDIF
\STATE Choose step-size $\alpha$, e.g. using \eqref{eq:armijo_lift}
\STATE $x_{k+1} \gets \Ret_{x_k}( - \alpha s_k)$
\ENDFOR
\end{algorithmic}
\end{algorithm}

\Cref{alg:rgd} provides pseudo-code for RGD on TTNs. 
The projector P$_x$ can be chosen $\Proj_x$ to iterate according
to \eqref{eq:rgd_H*} or $\Proj^\cartesian_x$ for \eqref{eq:rgd_H}. 
In the second case, we do not need to project $\nabla \bar f$
onto $\mathcal{T}$ first, 
because $\Proj^\cartesian_x$ can be extended smoothly to all of $\mathcal E$, where it holds that 
$\Proj^\cartesian_x = \Proj^\cartesian_x \circ \Proj_x$. Finally, we might skip
the projection all together, taking 
P$_x$ = Id, to iterate using the Euclidean gradient. In this case we still retract to keep all 
iterates orthogonal, avoiding numerical issues.

\subsection{Riemannian Newton's method \& trust-region}
Given an initial point $[x_0] \in \mathcal{T/G}$ and a retraction $\hat \Ret$,
RNM iteratively computes steps $[x_{k+1}] = \hat \Ret_{[x_k]}(\xi_{[x_k]})$ satisfying the Newton equation
\begin{equation} \label{eq:newton_equation}
    \Hess \hat{f}_{[x]}(\xi_{[x]}) = - \grad \hat{f}([x]).
\end{equation}
Rather than solving this equation directly in the quotient, we formulate its lift
to the horizontal space $H_x\mathcal{T}$:
\begin{equation} \label{eq:newton_equation_lift}
\begin{split}
    \Proj^H_x(\mathrm H_x(\xi_x)) \overset{\eqref{eq:hessian_lift}}{=} \ \ & \lift_x(\Hess \hat{f}_{[x]} (\xi_{[x]})) \\ 
    = - &\lift_x(\mathrm{grad} \hat{f}([x])) \overset{\eqref{eq:gradient_lift}}{=} -\Proj_x^H(\mathrm{grad} f(x)).
\end{split}
\end{equation} 
If \eqref{eq:newton_equation} has a unique solution $\xi_{[x]}$, then \eqref{eq:newton_equation_lift} also has 
a unique solution $\xi_x = \lift_x(\xi_{[x]})$ in $H_x\mathcal{T}$.
Assuming $\Hess \hat{f}_x$ is invertible at all iterates and that the retraction is equivariant, 
\cite[Thm. 6.3.2]{Absil:2008} guarantees at least quadratic convergence of 
the RNM in a neighborhood around a critical point, 
independently of the connection used to induce the covariant Hessian. This discussion does not 
only apply to $\HSt$, with $\mathrm{H}_x^\cartesian$ and $\Proj_x^\cartesian$, but also to $\HT$.
When the orthogonal horizontal space is employed, \eqref{eq:newton_equation_lift} simplifies to
\begin{equation} \label{eq:newton_equation_lift*}
    \Proj^{\times}_x(\Hess f_x(\xi_x)) = - \grad f(x).
\end{equation}
Although the projection on the RHS may be omitted, 
this formula still involves the costly evaluation of $\Proj^{\times}_x$ on the LHS,
when the lifted Newton equation is solved.
On the upside, \cite[Thm. 6.3.2]{Absil:2008} provides better convergence bounds, when the Riemannian connection is used.
Another alternative consists in simply solving 
the Newton equation of the total space
\begin{equation} \label{eq:newton_equation_total}
    \Hess f_x(\xi_x) = -\grad f(x).
\end{equation}
Existence and uniqueness of a solution is unclear for this system, and $\Hess f$ 
converges to a singular map, when approaching critical points \cite[Lemma 9.41]{Boumal:2023},
because the vertical space $\VT$ will be included in its kernel.
Nevertheless, one can find an approximate solution to \eqref{eq:newton_equation_lift*}
by running CG on \eqref{eq:newton_equation_total} and returning the previous CG-iterate 
upon failure. This strategy works because at critical points, CG initially spans
Krylov subspaces in $\HT$ for the first $\dim(\HT)$ iterations, before breaking as soon
as the kernel $\VT$ is reached \cite[Exercise 9.48]{Boumal:2023}.

\begin{algorithm}  
\caption{Riemannian trust-region}
\label{alg:rnm}
\begin{algorithmic}
\REQUIRE Objective function $f: \mathcal{T} \to \R$,
initial TTN-Parameter $x_0 \in \mathcal{T}$, tolerance $\varepsilon$, \\
pair of functions $(\mathrm{P}_x, \mathrm{H}_x)$,
retraction R$_x$
\FOR{$k$ = $0, 1, 2, \ldots$}
\STATE $\nabla \Bar{f} \gets $EuclideanGradient$(f, x_k)$
\STATE $s_k \gets  \mathrm{P}_{x_k}(-\nabla \Bar{f})$
\IF{$||s_k||^2 < \varepsilon$}
\RETURN $x_k$
\ENDIF
\STATE Approximately solve $\mathrm{H}_{x_k} d_k = s_k$ for $d_k$ using Steihaug--Toint-CG
\IF{$d_k$ meets trust-region criteria}
\STATE $x_{k+1} \gets \Ret_{x_k}(d_k)$
\ELSE 
\STATE $x_{k+1} \gets x_k$
\ENDIF
\ENDFOR

\end{algorithmic}
    
\end{algorithm}

Combining all three approaches gives rise to \Cref{alg:rnm}. 
The tuple of functions $(\mathrm{P}_x, \mathrm{H}_x)$ can be taken as $(\Proj^\cartesian_x, \Proj^\cartesian_x\circ\mathrm{H}^\cartesian_x)$
to accomplish \eqref{eq:newton_equation_lift}, $(\Proj_x, \Proj^\times_x\circ\Hess f_x)$ to realize \eqref{eq:newton_equation_lift*} or simply
$(\Proj_x, \Hess f_x)$ to solve \eqref{eq:newton_equation_total} on the total space $\mathcal{T}$.
Due to the inability of Newton's method to handle regions with indefinite Hessians, 
we do not use RNM in our implementation; 
instead, we employ a Riemannian trust-region algorithm, 
with Steihaug--Toint-CG as a subproblem solver \cite[Sec. 6.4]{Boumal:2023}. 
It can be seen as an extension to RNM, which only approximately solves the Newton equation in each iteration. At the same time 
it keeps track of the validity of the local quadratic model by restricting the update step to lie within a trust region, 
whose size is adjusted based on how well the model predicts the actual objective function decrease.
The trust-region mechanism was intentionally described in a simplified manner to improve readability. 
Note that both $\Hess f_x$ and $\Proj^\cartesian_x\circ\mathrm{H}^\cartesian_x$ serve as practical heuristics when employed as trust-region models. 
While convergence theory for RNM applies to $\Proj^\cartesian_x\circ\mathrm{H}^\cartesian_x$, 
theory for RTR only covers symmetric systems \cite[Sec. 7]{Absil:2008}.

\subsection{Comparison to optimization on low-rank manifolds}
In this subsection we want to draw up a brief comparison to Riemannian optimization on low-rank manifold $\mathcal{H}$.
As seen in \cite[Sec. 3.3]{Steinlechner:2016}, standard Riemannian gradient descent for objective function
$h : \R^{n_1 \times \ldots \times n_d \times K} \to \R$ would employ directions $\xi \in \HSt$ that accommodate
\begin{equation*}
    \Proj^\mathcal{H}_\mathbf{X}(\nabla h(\phi(x))) = \dd \phi(x)[\xi],
\end{equation*}
where $\Proj^\mathcal{H}_\mathbf{X}: \R^{n_1 \times \ldots \times n_d \times K} \to T_\mathbf{X}\mathcal{H}$ is the orthogonal
projector onto the tangent space of $\mathcal{H}$ at $\mathbf{X} = \phi(x)$. 
Update steps are then taken according
to $\dd \phi(x)[\xi]$, which has an efficient representation in low-rank format, 
requiring at most twice the bond dimensions
of the original network. A retraction on $\mathcal{H}$ is given by direct addition and HOSVD truncation \cite[Alg. 6]{Kressner:2012},
which represents a considerable computational effort. Furthermore, we want to point out that to our best knowledge, 
the details of optimization on $\mathcal{H}$ have only ever been fleshed out for tensor trains but not for TTNs. 

Our methods, on the other hand, iterate directly in $\mathcal{T}$ according to horizontal directions $\xi = \Proj^H_x(\dd \phi(x)^*[\nabla h(\phi(x))])$,
and retractions come at a comparably cheap cost. 
Additionally, we want to report on a peculiarity involving the horizontal space $\HSt$. 
Recall from Section \ref{sec:orthogonal_ttns} that the bond dimensions of full-rank TTN parameters $\mathcal E^*$ satisfy
$k_t \leq k_{\tl} k_{\tr}$.
The case $k_t = k_{\tl} k_{\tr}$ is a practical choice for bond dimensions, especially on the lower levels of the TTN. 
Now, we have previously mentioned, that $\HSt$ reads as a Cartesian product of horizontal spaces $H_{B_t}\St(k_{\tl}k_{\tr}, k_t)$. For the transfer tensors 
$B_t \in \R^{k_{\tl}k_{\tr} \times k_t}$ that fulfill the above condition, 
according to \cite[Sec. 9.16]{Boumal:2023} these horizontal spaces are of dimension 
\begin{equation} \label{eq:pooling_layers}
    k_t(k_{\tl}k_{\tr} - k_t) = 0.
\end{equation}
Quotient optimization involving the Cartesian horizontal space will not modify those tensors, 
and update steps on them may be omitted entirely. This leads to a surprising conclusion: Due to the gauge freedom
of TTNs, those tensors need not be optimized at all, if the other tensors are suitably optimized. While the horizontal space
$\HSt$ successfully leverages this fact, $\HT$ intrinsically couples all transfer tensors, and fails to isolate unneeded degrees 
of freedom to individual tensors. Even though Riemannian optimization on $\mathcal{H}$ theoretically offers the same advantage, 
a computational trick that is used in practice to redeem numerical issues \cite[Sec. 3.3]{Steinlechner:2016} reintroduces coupling
through repeated reorthogonalizations. The original update direction $\xi$ is shifted to a different horizontal space, which still 
resembles but is not equivalent to $\HSt$, and the described effect cannot be observed. 

Finally, a comment about the closure of $\mathcal{T}$ is in order. 
Since $\mathcal{T}$ forms as an intersection with full-rank parameters $\mathcal{E}^*$,
it is not closed in $\mathcal{E}$, and the accumulation points of series in $\mathcal{T}$ might not
be full-rank themselves, thus breaking convergence results that were mentioned earlier.
For more details on the closure problem, we refer to 
\cite[Sec. 4.6]{Uschmajew:2013}.
In practice, this issue can be addressed by descending to a manifold
 of
lower HT-rank through the means of HOSVD truncation \cite[Alg. 6]{Kressner:2012},
whenever 
HT singular values numerically approach zero, and proceeding the iteration from there.
It should be noted, however, that simple rank-adaptive schemes based on singular value thresholds may not be sufficient to 
prevent overfitting in supervised learning tasks. 
A more robust dynamical rank adaptation remains an open area of research in tensor network learning.

\section{Machine Learning with TTNs} \label{sec:ml_ttn}
In this section we explore how Riemannian optimization on TTNs 
can be employed for supervised machine learning. Following the seminal work of \cite{Stoudenmire:2016},
an input vector $v \in \R^d$ is encoded as a tensor product
\begin{equation} \label{eq:product_state}
    \Phi(v) = \varphi^1(v_1) \otimes \ldots \otimes \varphi^d(v_d),
\end{equation}
where $(\varphi^i)_{i  \in \{1, \ldots, d\}}, \varphi^i: \R \to \R^{n_i}$ 
are local feature maps applied to each input component separately.
The model is parametrized by a $d+1$-dimensional tensor $\mathbf{X} \in \R^{n_1 \times \ldots \times n_d \times K}$, and the model equation reads
\begin{equation*}
    y^j(v) =  \dotpro{\varphi^1(v_1) \otimes \ldots \otimes \varphi^d(v_d) \otimes e_j}{\mathbf{X}}.
\end{equation*}
This model is intended to provide a classifier with $K$ classes,
and a classification for input $v$ may be evaluated as $\arg \max_{1 \leq j \leq K}(y^j(v))$ \cite[Eq. 4]{Stoudenmire:2016}.
In matricized form, the above equation gives rise to the \emph{response vector} $y = (y^j(v))_{1 \leq j \leq K} \in  \R^{K}$,
which contains all model responses for a given input vector:
\begin{equation*} 
    y^\T = (\varphi^1(v_1)^\T \otimes \ldots \otimes \varphi^d(v_d)^\T \otimes I)\mathbf{X} = \vect{\otimes_{1 \leq i \leq d} \varphi^i(v_i)}^\T \mathbf{X}^{(1, \ldots, d)}.
\end{equation*}
Memory requirements for both the product state \eqref{eq:product_state} and tensor $\mathbf{X}$ are prohibitive, which is
why we represent $\mathbf{X}$ by a tree tensor network. 
Specifically, we employ an orthogonal TTN-parameter $x \in \mathcal{T}$ to represent $\mathbf{X} = \phi(x)$. The feature maps 
$\varphi^i: \R \to \R^{n_i}$ determine the external dimensions $n_t = n_i$ for $t = \{i\}$, and the number of classes used determines the label dimension
$k_{\tb} = K$ of the TTN. Importantly, the bond dimensions, $(k_t)_{t \in J^-}$ govern the computational complexity of the model. 
They must be sufficiently large to represent the underlying data accurately, yet small enough to keep the computational cost manageable. 
Additionally, regularization through small bond dimensions can help prevent overfitting, 
although tensor networks have been observed to exhibit a second descent in the bias-variance tradeoff
for vastly overparametrized models, similar to neural network models \cite{Pomarico:2025}.
For our application, a suitable strategy is given by defining some maximum bond dimension $k_{t_{\mathrm{max}}}$ and setting $k_t = \min\{k_{\tl}k_{\tr}, k_{t_{\mathrm{max}}}\}$,
or by using Stoudenmire's unsupervised construction algorithm \cite{Stoudenmire:2018}. 
Note that both approaches lead to an abundance of transfer tensors that realize $k_t = k_{\tl}k_{\tr}$, subjugating them to the effects described 
below \Cref{eq:pooling_layers}.

For evaluating a model response or performing optimization routines on the model, 
neither the high-dimensional tensor $\mathbf{X}$ nor the tensor product \eqref{eq:product_state} have to be formed explicitly. Instead we will have to consider
vectors and tensors along the TTN, which depend on the bond dimensions.
\subsection{Forward propagation}
From now on, we assume that any input data $v$ has already been prepared by feature maps $\varphi^1, \ldots, \varphi^d$, forming
an \emph{input sample} $\mathbf{s} = (s_i)_{i \in D}, s_i \in \R^{n_i}$ with $s_i = \varphi^i(v_i)$. Given an input sample and
a TTN associated with tensor $\mathbf{X}\in\R^{n_1 \times \ldots \times n_d\times K}$, the response vector is calculated as
\begin{equation} \label{eq:contract_sample}
    y^\T =
    \vect{\otimes_{i \in \tb} s_i}^\T U_{\tb}.
\end{equation}
By induction on \eqref{eq:ttn_recursion}, the response vector can be computed recursively. Starting from the leaves, 
calculate \emph{effective samples} $(s_t)_{t \in J}, s_t \in \R^{k_t}$ with
\begin{align} \label{eq:forward_prop}
\begin{split}
    &s_t^\T = \vect{\otimes_{i \in t} s_i}^\T U_t  = 
    (s_{\tl}^\T \otimes s_{\tr}^\T \otimes I) \mathbf{B}_t \text{ for } t \in J,
\end{split}
\end{align}
yielding $y = s_{\tb}$ at the root node. This formula represents a generalization from e.g. \cite[Thm. 1]{Novikov:2017}
to tree tensor networks, and
we call \eqref{eq:forward_prop} the \emph{forward propagation} of sample $\mathbf{s}$.

\subsection{Backpropagation}
Given an \emph{expected response} $y^* \in \R^K$ associated with input sample $\mathbf{s}$, 
define a loss function 
$\mathcal{L}: \R^K \times \R^K \to \R, (y, y^*) \mapsto \mathcal{L}(y, y^*)$. 
For brevity, we drop the explicit dependence of $\mathcal{L}$ on $y^*$ 
and assume $y^*$ is fixed. Plugging \eqref{eq:contract_sample} into $\mathcal{L}$ then yields the function 
\begin{equation} \label{eq:objective_function_ml}
    h: \R^{n_1 \times \ldots \times n_d\times K} \to \R, \qquad \mathbf{X} \mapsto \mathcal{L}((s_1^\T \otimes \ldots \otimes s_d^\T \otimes I)\mathbf{X}).
\end{equation}
An example for a loss function would be the \emph{quadratic error loss function} 
$\mathcal{L}_2(y, y^*) = \tfrac{1}{2}||y^* - y||^2$ , 
for which the gradient simply reads $\nabla \mathcal{L}_2(y) = y - y^*$. 

Returning to the optimization problem \eqref{eq:optimization_problem}, 
our goal is to minimize $\hat{f} = h \circ \hat{\phi}$ on $\mathcal{T/G}$. 
To find an expression for $\grad \hat f$, we derive the Euclidean gradient $\nabla \Bar f$ for our specific choice of $h$.
The Euclidean gradient of $h$ at $\mathbf{X}$ reads
\begin{equation*}
    \begin{split} 
    \nabla h(\mathbf{X}) = s_1 \otimes \ldots \otimes s_d \otimes \nabla \mathcal{L}(y) =  s_1 \otimes \ldots \otimes s_d \otimes \ell,
\end{split}
\end{equation*}
where we defined the \emph{loss gradient} $\ell = \nabla \mathcal{L}(y)$, 
and employed the multivariate chain rule:
$\nabla (k \circ g)(x) = g'(x)^\T \nabla k(g(x))$ for $g: \R^n \to \R^m, k: \R^m \to \R$.
Additionally define \emph{effective loss gradients} $(\ell_t)_{t \in J}, \ell_t \in \R^{k_t}$ given by
    \begin{align}\label{eq:backward_prop}
\begin{split}
    &\ell_{\tl} = ( I \otimes s_{\tr}^\T \otimes \ell_t^\T ) \mathbf{B}_t,  \\
    &\ell_{\tr} =  (  s_{\tl}^\T \otimes I \otimes \ell_t^\T ) \mathbf{B}_t,
\end{split}
\end{align}
which can be calculated recursively from the loss gradient $\ell = \ell_{\tb}$ at the root node.
We call \eqref{eq:backward_prop} the \emph{backpropagation} of loss gradient $\ell$, 
because it allows to construct the Euclidean gradient of $f$ through what is 
really a repeated application
of the chain rule.

\begin{figure}[h] 
\centering
\includegraphics[width=0.88\textwidth]{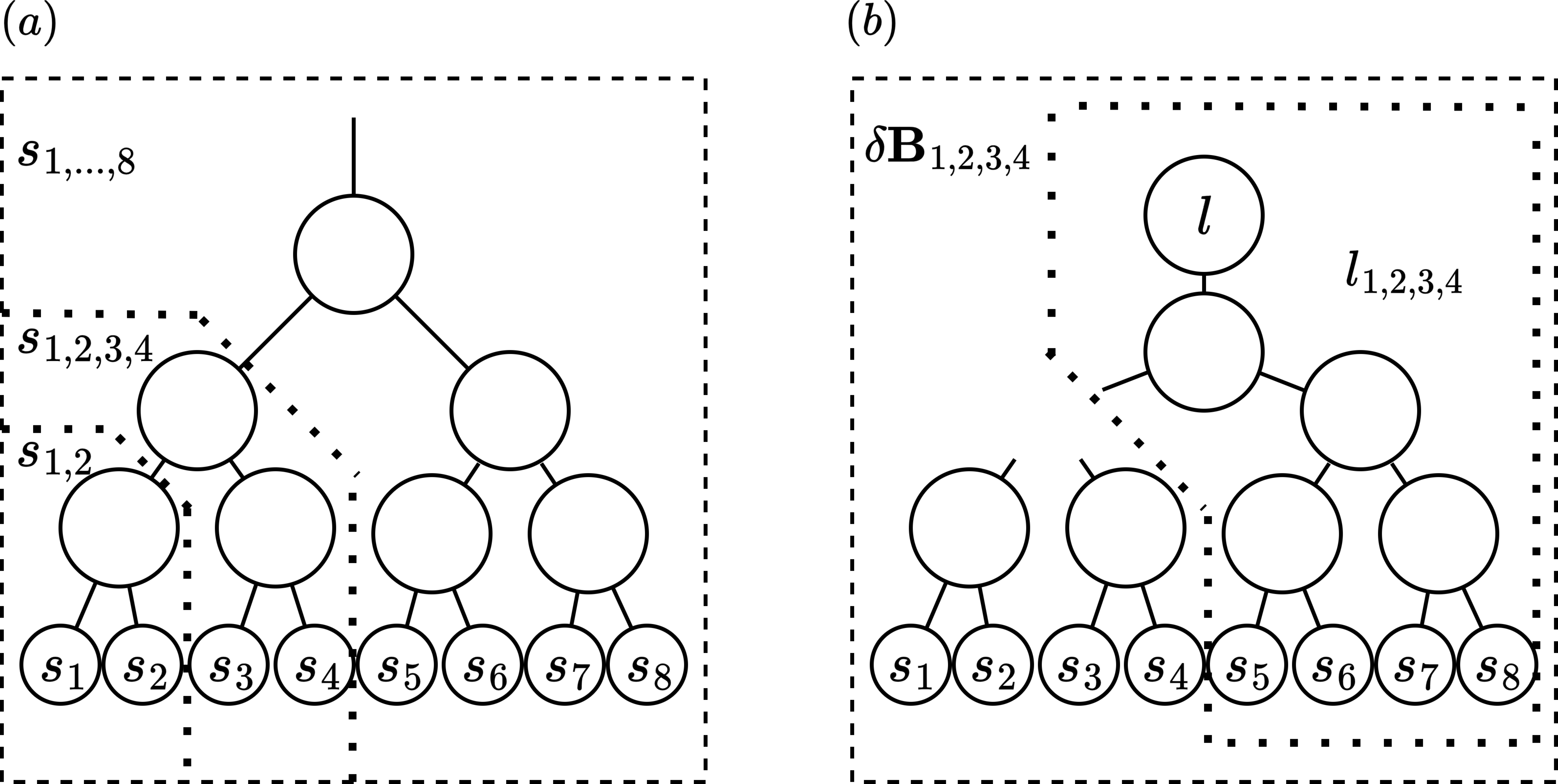} 
\caption{forward propagation of a sample and component $\delta \mathbf{B}_{1,2,3,4}$ of the Euclidean gradient}
\label{fig:ml_ttn}
\end{figure}

\begin{theorem}
    Let $x = (\mathbf{B}_t)\in \mathcal{T}$ and let $\ell_t$ be the backpropagated gradients \eqref{eq:backward_prop} of function $h$ as defined in \eqref{eq:objective_function_ml}.
    Then $\nabla \bar f(x)  = \nabla (h \circ \phi)(x) = \dd \phi(x)^*[\nabla h(\phi(x))]$ reads
    \begin{equation} \label{eq:gradient_outer_product}
        \nabla \bar f(x) = (\delta \mathbf{B}_t) = (s_{\tl} \otimes s_{\tr} \otimes \ell_t)_{t \in J}.
    \end{equation}
\end{theorem}
\begin{proof}
    Let $\nabla \bar f(x) = (\delta \mathbf{B}_t)$ be accompanied by tensors $\mathbf{W}_t$ 
    like in \eqref{eq:gradient_recursion}. Assume that for some $t \in J$, it has the form 
    $\mathbf{W}_t = \vect{\otimes_{i \in \tl}s_i} \otimes \vect{\otimes_{i \in \tr}s_i} \otimes \ell_t$.
    Now calculate
    \begin{align*} 
\begin{split}
    \delta \mathbf{B}_t &\overset{\eqref{eq:gradient_recursion}}{=} (U_{\tl}^\T \otimes U_{\tr}^\T \otimes I) \mathbf{W}_t
    = U_{\tl}^\T \vect{\otimes_{i \in \tl}s_i} \otimes U_{\tr}^\T \vect{\otimes_{i \in \tr}s_i} \otimes \ell_t \\
    &\overset{\eqref{eq:forward_prop}}{=} s_{\tl} \otimes s_{\tr} \otimes \ell_t ,
    \\
    \delta U_{\tl} &\overset{\eqref{eq:gradient_recursion}}{=} \dotpro{U_{\tr}^\T \times_2 \mathbf{W}_t}{\mathbf{B}_t}_{(2, 3)}
    \overset{[9]}{=} \vect{\otimes_{i \in \tl}s_i} \otimes \left( \vect{\otimes_{i \in \tr}s_i}^\T U_{\tr} \times_2 \ell_t^\T \times_3 \mathbf{B}_t \right) \\
    &\overset{\eqref{eq:forward_prop}}{=} \vect{\otimes_{i \in \tl}s_i} \otimes \left( I\otimes s_{\tr}^\T \otimes \ell_t^\T \right) \mathbf{B}_t 
    \overset{\eqref{eq:backward_prop}}{=} \vect{\otimes_{i \in \tl}s_i} \otimes \ell_{\tl},
    \\
    \delta U_{\tr} &\overset{\eqref{eq:gradient_recursion}}{=} \dotpro{U_{\tl}^\T \times_1\mathbf{W}_t}{\mathbf{B}_t}_{(1, 3)}
    \overset{[9]}{=} \vect{\otimes_{i \in \tr}s_i} \otimes \left( \vect{\otimes_{i \in \tl}s_i}^\T U_{\tl} \times_1 \ell_t^\T \times_3 \mathbf{B}_t \right) \\
    &\overset{\eqref{eq:forward_prop}}{=} \vect{\otimes_{i \in \tr}s_i} \otimes \left(  s_{\tl}^\T \otimes I \otimes \ell_t^\T \right) \mathbf{B}_t 
    \overset{\eqref{eq:backward_prop}}{=} \vect{\otimes_{i \in \tr}s_i} \otimes \ell_{\tr}.
\end{split}
\end{align*}
In the fifth and ninth equation we used \cite[Eq. 2]{DaSilva:2015}, which states that multilinear factors, whose mode is
not involved in the contraction can be factored out. 
Finally recall from \eqref{eq:gradient_recursion}, 
that $W_{\tb} = (\nabla h (\phi(x)))^{(D)} = \vect{\otimes_{i \in \tb}s_i} \otimes \ell_{\tb}$,
then the assertion holds by induction.
\end{proof}
Combining forward propagation \eqref{eq:forward_prop} and backward propagation \eqref{eq:backward_prop}
constitutes Algorithm \ref{alg:ttn_backprop}, which computes $\nabla \bar f$ at any $x \in \mathcal{T}$. This algorithm
also has a nice interpretation in tensor diagram notation. \Cref{fig:ml_ttn}(a) represents the forward pass, 
allowing to calculate $\ell = \nabla\mathcal{L}(s_{\tb})$ at the top. Then the gradient \eqref{eq:gradient_outer_product} for any node is given by removing that node from 
the network, and contracting everything else, as exemplified in \Cref{fig:ml_ttn}(b).
\begin{algorithm}[H]

\caption{Euclidean gradient for TTNs}
\label{alg:ttn_backprop}
\begin{algorithmic}
\REQUIRE TTN-Parameters $x = (\mathbf{B}_t) \in \mathcal{T}$, input sample $\mathbf{s}$, loss function $\mathcal{L}$
\FOR{$t \in J$, visiting children before their parents}
\STATE $s_t \gets (s_{\tl}^\T \otimes s_{\tr}^\T \otimes I) \mathbf{B}_t$
\ENDFOR
\STATE $\ell_{\tb} \gets \nabla\mathcal{L}(s_{\tb})$
\FOR{$t \in J^-$, visiting parents before their children}
\STATE $\ell_{\tl} \gets \left( I \otimes s_{\tr}^\T \otimes \ell_t^\T \right) \mathbf{B}_t$
\STATE $\ell_{\tr} \gets  \left(  s_{\tl}^\T \otimes I \otimes \ell_t^\T \right) \mathbf{B}_t$
\ENDFOR
\FOR{$t \in J$}
\STATE $\delta \mathbf{B}_t \gets s_{\tl} \otimes s_{\tr} \otimes \ell_t $
\ENDFOR
\RETURN  $(\delta \mathbf{B}_t)$
    
\end{algorithmic}
    
\end{algorithm}

\subsection{Implementation details}

For any real machine learning task, 
we will have multiple pairs $(y^*_j, \mathbf{s}^j)_{1 \leq j \leq m}$ of samples and associated objectives. 
Let $y_j$ denote the response generated by sample $\mathbf{s}^j$ via \eqref{eq:contract_sample}. 
The (total) loss function then takes the form
\begin{equation*}
    \mathcal{L}^m: (\R^K \times \R^K)^m \to \R, (y_j, y^*_j)_{1 \leq j \leq m} \mapsto \frac{1}{m} \sum_{j=1}^m \mathcal{L}(y_j, y^*_j)
\end{equation*}
yielding a function $h^m: \R^{n_1 \times \ldots \times n_d\times K} \to \R$
\begin{equation*}
    h^m(\mathbf{X}) = \frac{1}{m} \sum_{j=1}^m \mathcal{L}((s_{1}^j \otimes \ldots \otimes s_{d}^j \otimes I)^\T\mathbf{X}, y^*_j),
\end{equation*}
which we can use for optimization on our TTN. Due to the linearity of gradients, 
we can evaluate both $h^m$ and $\nabla \bar f$ at any point 
separately for all samples and accumulate afterwards. 

This means that both forward and backward propagation require evaluating quantities of the form $s_t^j = (s_{\tl}^j \otimes_K s_{\tr}^j)^\T B_t$,
where $B_t$ remains fixed, and $(s_{\tl}^j$, $s_{\tr}^j)_{1 \leq j \leq m}$ represent pairs of vectors, that need to be contracted. 
Arranging the individual vectors columnwise into matrices
$S_\tl = (s_{\tl}^1, \ldots,  s_{\tl}^m) \in \R^{k_\tl \times m}$, 
$S_\tr = (s_{\tr}^1, \ldots, s_{\tr}^m) \in \R^{k_\tr \times m}$ 
and $S_t = (s_{t}^1, \ldots, s_{t}^m) \in \R^{k_t \times m}$ allows to formulate the problem as 
\begin{equation*}
    S_t^\T = (S_\tl \odot S_\tr)^\T B_t,
\end{equation*}
which is the well-known matricized-tensor-times-Khatri-Rao-product (MTTKRP) \cite{Bader:08}, and which we implemented in a parallelized fashion
using BLAS routines. A similar strategy applies to forming and summing over the outer products of \eqref{eq:gradient_outer_product}. It holds that
\begin{equation*}
    \delta B_t = \frac{1}{m}\sum_{j=1}^m (s_{\tl}^j \otimes_K s_{\tr}^j) \otimes \ell_t^j = \frac{1}{m}(S_\tl \odot S_\tr) L_t^\T,
\end{equation*}
where the backpropagated loss gradients $\ell_t^j$ form the columns of $L_t = (\ell_t^1, \ldots, \ell_t^m) \in \R^{k_t \times m}$. Again this operation is implemented in parallel using
BLAS routines. 
Both of the above procedures represent the workhorses, that drive our algorithms. In order to achieve suitable performance,
they are implemented in a dedicated \texttt{C++} library. This also includes QR, SVD and matrix inversion algorithms, that are required by our retractions, and which are implemented through LAPACKE.

\Cref{alg:rgd,alg:rnm,alg:cartesian_retraction} themselves are implemented in $\texttt{python}$, 
delegating only performance critical operations to \texttt{C++}. 
This includes the projectors $\Proj, \Proj^\equiv$ and $\Proj^\times$ \eqref{eq:proj_H*},
where the latter requires esolving the pseudoinverse of $\Proj^{H,V}_x \circ (\Proj^{H,V}_x)^*$ in $\HSt$. 
Since the operator is positive definite on $\HSt$,
our implementation uses the \texttt{numpy} CG solver for this task. 
The system is not solved exactly, but rather up to a relative decrease of residuals to $10^{-2}$.

Initial iterates $x \in \mathcal T$ for our experiments were found by randomly initializing $\bar x\in \mathcal E^*$ and applying the orthogonalization algorithm \cite[Alg. 1]{Kressner:2012}.
The line-search for first-order algorithms is implemented using two-way backtracking with Armijo parameter $r=10^{-4}$. Parameters and other practical advice for the trust-region mechanism and Steihaug--Toint-CG are implemented exactly as they
appear in \cite[Sec. 6.4.6]{Boumal:2023}. 
The operators $\mathrm{H}^\cartesian$ and $\Hess f$ are approximated using 
their respective gradients in
a finite difference approach \cite[Ex. 5.32]{Boumal:2023}.

\section{Numerical Evaluation}
\label{sec:results}

We compare the various optimization approaches by training a TTN on 
the \texttt{digits} dataset \cite{digits:1998}, comprising $N=1797$ 8$\times$8 grey scale images 
of handwritten digits. Consequently, the number of external sites is chosen $d = 64$. 
The objective function $f$ is given via the quadratic error loss $\mathcal{L}_2$.
We use the spin feature map $\varphi: x \mapsto (\cos(\tfrac{\pi}{2}x), \sin(\tfrac{\pi}{2}x))$ 
as proposed in \cite[Eq. 3]{Stoudenmire:2016}, applied to
each gray scale pixel value $v_i \in [0, 1]$ of an input vector $v$.
Thus the external dimensions must be $n_i = 2$ for all $i \in \{1, \ldots, d\}$. 
The bond dimensions are chosen $k_t = \min\{n_t, 8\}$ for all $t \in J^-$.
The label dimension is chosen $K = 10$, and the expected responses $y^*_n$ are taken one-hot vectors $e_{j+1}$, 
with $j \in \{0, \ldots 9\}$ the handwritten digit represented by the corresponding sample $\mathbf{s}_n$. 
The choice of bond dimensions is consistent with what previous work \cite[Sec. VI]{Stoudenmire:2016} 
suggests for tensor network learning of similarly difficult datasets, 
and provides a sufficient balance between model expressivity and computational efficiency.
This results in effectively $\dim(\mathcal{T/G})=7936$ degrees of 
freedom, which is on the same order as the number of parameters used in typical models for classification on the \texttt{digits} dataset. 
For instance, the standard multilayer perceptron configuration in \texttt{scikit-learn} uses a single
hidden layer of 100 neurons resulting in approximately 
$7400$ parameters \cite{scikit-learn}. 
We verified that the performance of presented optimization techniques is 
qualitatively consistent for other dimensions (e.g., $k_{t_\mathrm{max}} \in \{4,16\}$); 
therefore, $k_{t_\mathrm{max}}=8$ is used throughout this work as a 
representative baseline to evaluate the optimization algorithms.

As is common for machine learning tasks, only $80 \%$ of samples are used for optimization, 
forming the training set, the remainder is used for testing. 
All results are averaged over 10 runs, each defined by a pair of randomly generated initialization and train-test split. 
The same 10 pairs are used for all algorithms.
Throughout our evaluation we employed the Cartesian retractions induced by the QR decomposition (QR), 
the polar decomposition (PD) and the Cayley transform (CT) (see \Cref{thm:cartesian_retraction}).

\begin{figure}[h]
    \centering
    \includegraphics[trim={5cm 0 5cm 0}, width=\textwidth]{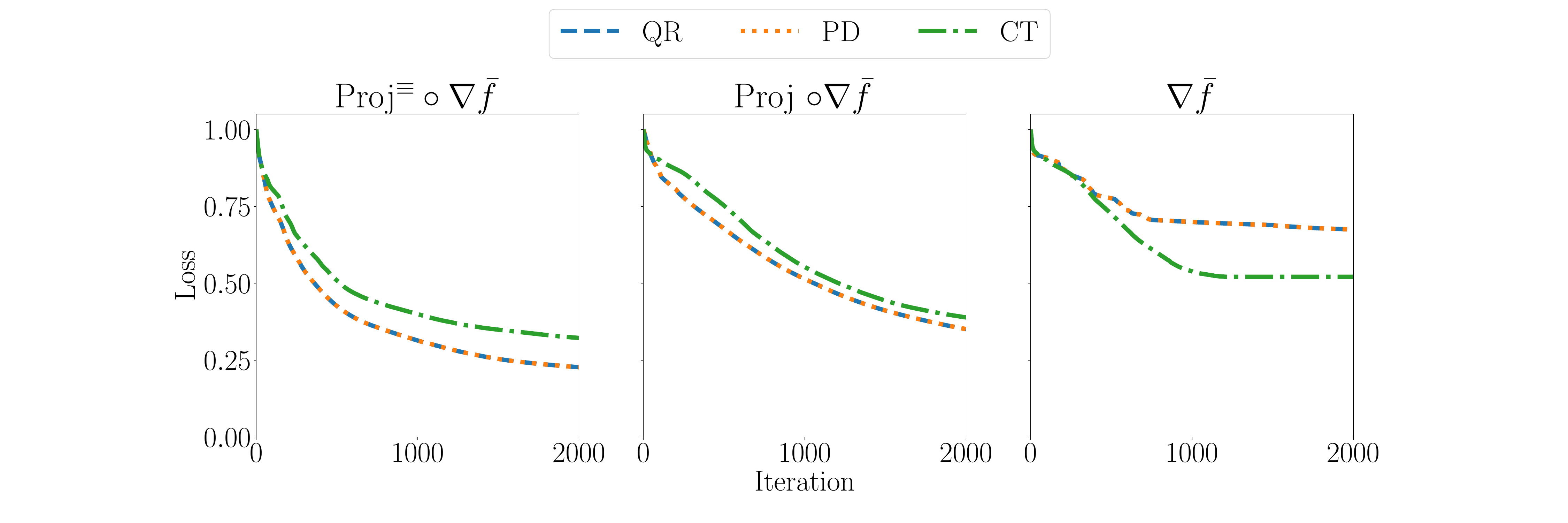}
    \caption{2000 iterations of RGD for diverse approaches and retractions}
    \label{fig:rgd}
\end{figure}

First, we compare the performance of RGD for $\Proj^\cartesian \circ \nabla\bar f = \lift^\cartesian \circ \grad^\cartesian \hat f$ \eqref{eq:rgd_H} to that of 
$\Proj \circ \nabla\bar f = \lift^\times \circ \grad^\times \hat f = \grad f$ \eqref{eq:rgd_H*}. We also include the Euclidean gradient $\nabla\bar f$ into our
benchmarks, which is achieved by skipping the projection in \Cref{alg:rgd}. As seen in \Cref{fig:rgd}, the choice of the retraction
plays a minor role for both Riemannian approaches. Even though the QR retraction does not formally satisfy the requirements of a quotient retraction, it remains
a viable numerical variant for quotient optimization.
The descent of both strategies is somewhat comparable, which is not surprising, 
since both of them are standard Riemannian gradient descent algorithms on the quotient, differing only in their metric, which does not directly influence
convergence guarantees \cite[Thm. 4.3.1]{Absil:2008}. Using the Euclidean gradient does not appear viable:
step-size control would usually get stuck after a while, because no step-size fulfilling the Armijo-condition could be found. 
This strongly indicates that the Euclidean
gradient does not necessarily deliver descent directions in combination with the retractions used, 
which in turn justifies usage of Riemannian optimization on TTN.

\begin{figure}[h]
    \centering
    \includegraphics[trim={5cm 0 5cm 0}, width=\textwidth]{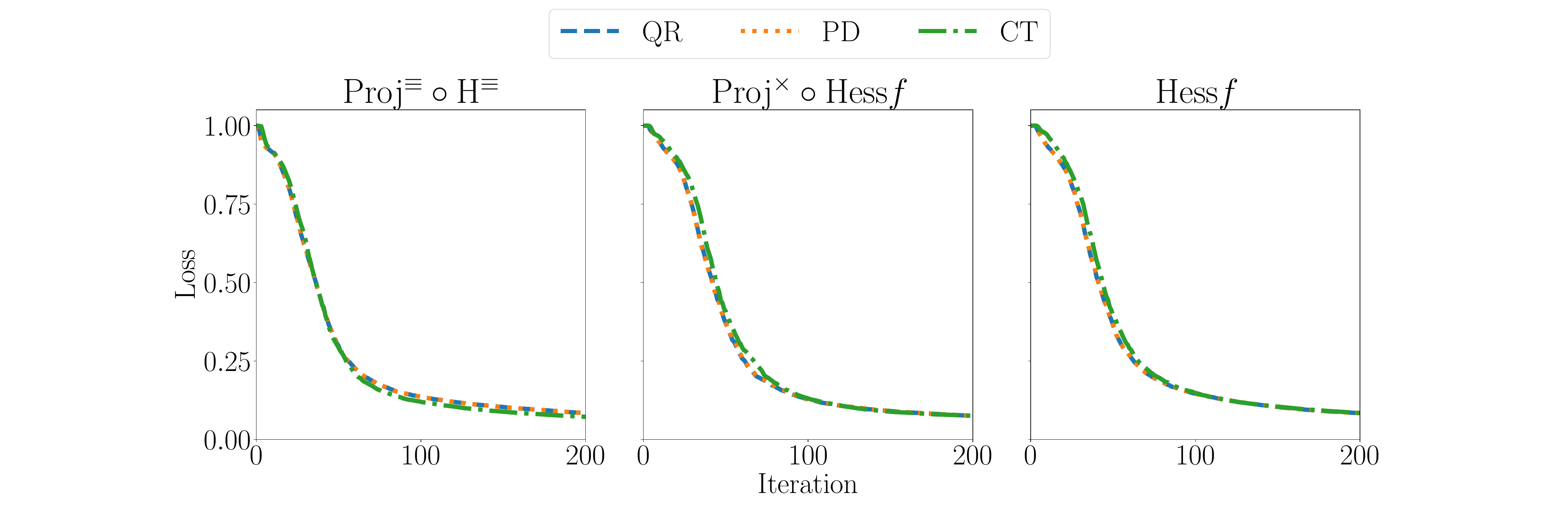}
    \caption{200 iterations of RTR for diverse approaches and retractions}
    \label{fig:rtr}
\end{figure}
We now move on to second-order optimization algorithms.
For the trust-region models of \Cref{alg:rnm}, we choose the quotient Hessians
$\Proj^\cartesian \circ \mathrm{H}^\cartesian$ \eqref{eq:newton_equation_lift}, 
$\Proj^{\times} \circ \Hess f$ \eqref{eq:newton_equation_lift*} and
$\Hess f$ \eqref{eq:newton_equation_total}, the Riemannian Hessian on $\mathcal{T}$. We exclude the Euclidean Hessian, as the corresponding algorithm did 
not converge at all in our experiments. Looking at \Cref{fig:rtr}, 
we observe that all three strategies show an almost identical descent for any retraction, 
again including the non-quotient QR retraction, and we report
that all approaches required a similar amount of inner CG iterations throughout their descent.
It was already asserted in \cite[Ex. 9.48]{Boumal:2023}, that $\Hess f$ tends to deliver results close
to those of $\Proj^\times \circ \Hess f$ in practical application, so this behavior is not entirely unexpected (see also the discussion below \eqref{eq:newton_equation_total}). Interestingly, 
$\Proj^\cartesian \circ \mathrm{H}^\cartesian$ performs on par with those approaches. We mentioned earlier that superlinear convergence 
of the quotient-based RNM can be guaranteed independently of the connection used, 
but as seen in the proof of \cite[Thm. 6.3.2]{Absil:2008}, using a Riemannian
connection offers better theoretical convergence bounds. Furthermore, even though RTR can formally only handle symmetric systems,
it proved effective for non-symmetric $\mathrm{H}^\cartesian$ with no further modification. 
A possible explanation might be, that $\mathrm{H}^\cartesian$ is close to symmetric, and Steihaug--Toint-CG handles indefinite models robustly anyways.
\begin{figure}[h]
    \centering
    \includegraphics[width=1\textwidth]{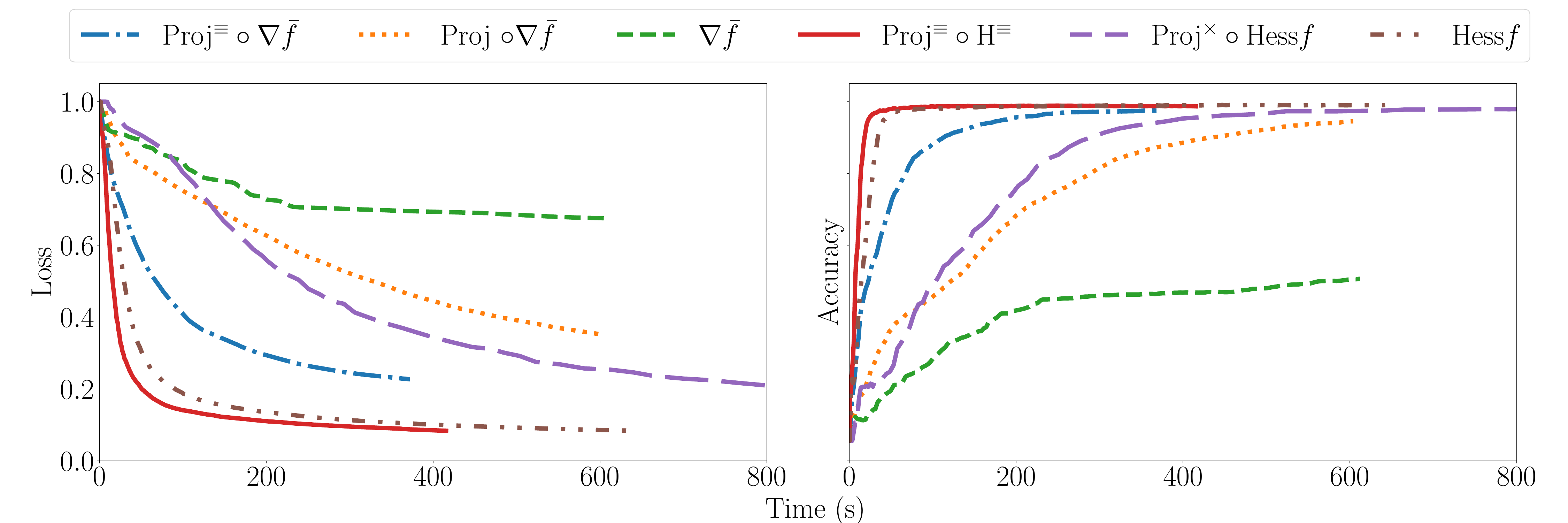}
    \caption{Performance comparison of RGD and RTR using QR retraction. 
    RGD/RTR algorithms were stopped after 2000/200 iterations respectively.}
    \label{fig:performance}
\end{figure}

It is to be expected, that RTR converges in fewer iterations than RGD, but the per-iteration complexity is much higher in the former. 
To ensure a fair comparison, \Cref{fig:performance} shows the progress of algorithms over their runtime. 
RTR algorithms easily outperform RGD, except for $\Proj^\times \circ \Hess f$, which requires the costly projection onto the orthogonal horizontal space.
Routines involving the horizontal space $\HSt$ have a significant advantage over their counterparts.
For example, $\Proj^\cartesian \circ \mathrm{H}^\cartesian$ completes 200 iterations in only 60\% of the time required by $\Hess f$.
This happens because working with $\HSt$ enables the exclusion of optimization updates on transfer tensors that fulfill \eqref{eq:pooling_layers}, 
resulting in numerically indistinguishable optimization iterates while significantly reducing computational load.
Analysis of classification accuracies over time suggests that the order of convergence 
to a local minimum does not necessarily improve the generalization capabilities of the model.
In the machine learning context the actual advantage of second-order methods 
is likely their faster initial convergence.

Finally, to validate the applicability of our algorithms for machine learning on TTNs, we provide \Cref{tab:accuracies},
that displays the percentage of correctly classified inputs from the test data, after training of the TTN was completed. 
All Riemannian approaches involving QR and PD rectractions yield accuracies exceeding 90\%, demonstrating the success of those training procedures.
The accuracies are comparable to those achieved by support vector machines and multilayer perceptrons
(97\% -- 99\% with \texttt{scikit-learn}), 
confirming that our model is competitive in terms of expressive power.
\begin{table} 
\centering
\begin{tabular}{c |c | c | c }
Method & QR & PD & CT \\
\hline
Proj$^\equiv \circ \nabla \bar f$ & $97.51 \pm 8.86 \times 10^{-3}\%$ & $97.51 \pm 8.86 \times 10^{-3}\%$ & $83.75 \pm 4.71 \times 10^{-0}\%$ \\
Proj $\circ \nabla \bar f$ & $94.52 \pm 2.97 \times 10^{-2}\%$ & $94.52 \pm 2.97 \times 10^{-2}\%$ & $83.29 \pm 8.10 \times 10^{-0}\%$ \\
$\nabla \bar f$ & $50.67 \pm 6.44 \times 10^{-0}\%$ & $50.67 \pm 6.44 \times 10^{-0}\%$ & $77.89 \pm 5.44 \times 10^{-0}\%$ \\
\hline
Proj$^\equiv \circ \mathrm H^\equiv$ & $98.69 \pm 2.17 \times 10^{-3}\%$ & $98.69 \pm 2.17 \times 10^{-3}\%$ & $98.89 \pm 1.59 \times 10^{-3}\%$ \\
Proj$^\times \circ \mathrm{Hess}f$ & $98.91 \pm 1.73 \times 10^{-3}\%$ & $98.97 \pm 1.65 \times 10^{-3}\%$ & $98.84 \pm 1.98 \times 10^{-3}\%$ \\
$\mathrm{Hess}f$ & $99.05 \pm 8.00 \times 10^{-4}\%$ & $99.05 \pm 8.00 \times 10^{-4}\%$ & $99.02 \pm 1.82 \times 10^{-3}\%$ \\
\end{tabular}
\caption{Mean and standard deviation of classification accuracies on test data after training}
\label{tab:accuracies}
\end{table}

\section{Conclusion and Outlook}

We devised the application of first- and second-order optimization algorithms for machine learning 
with tree tensor networks.

Arising from the non-uniqueness of TTN-parameters, 
we identified the quotient manifold of tree tensor networks as our optimization space of interest. 
We employed the orthogonal TTN-parameters, and vectors from two different horizontal distributions to serve as proxies
in the quotient optimization: the \emph{orthogonal} horizontal distribution $H^\times$ and the \emph{Cartesian} horizontal distribution $H^\equiv$. Explicit forms of important projectors adhering to those sets were presented. 
We demonstrated that these two distributions correspond to different metrics on the quotient space, 
which in turn induce different Riemannian gradients, connections and covariant Hessians on the quotient manifolds. 
Moreover, we conceived several efficient retractions for the manifolds of TTNs, which, assembled together with the other optimization
ingredients, resulted in multiple versions of Riemannian gradient descent and Riemannian trust-region algorithms. 
While the orthogonal horizontal distribution yields more concise formulas for those, 
the Cartesian horizontal distribution constitutes a viable alternative.
Finally, we provided a description of non-linear kernel learning with TTNs
and devised a backpropagation algorithm in this context. We numerically evaluated 
the presented algorithms on an image classification task. The results highlighted the importance of considering optimization on TTNs 
as Riemannian, instead of Euclidean, and displayed the strong advantage of including second-order information. Surprisingly, 
we could not observe any significant drawbacks in the practical application of the Cartesian horizontal space. It instead
gave rise to the best-performing algorithms.

A natural extension of our work would be to generalize it to more complex tensor network architectures, such as MERA, 
which have unique capabilities when used in machine learning \cite{Reyes:2020}. 
It would be interesting to see a comparison with
modern machine learning classifiers (such as convolutional neural networks) and the algorithms of the present work. 
To further enhance scalability and performance, future work could incorporate stochastic versions of 
Riemannian gradient descent \cite{Bonnabel:2013} and Riemannian trust-region methods.
Additionally, tensor decompositions have already been successfully used 
to compress trained neural networks \cite{ONeill:2020}, 
while offering a much richer model interpretability \cite{Felser:2021}. 
Therefore it would make
sense to combine both TTN compression techniques and explainability techniques 
with our training procedures.

\section*{Acknowledgements}
The authors would like to thank Timo Felser and Tensor AI Solutions for supporting this research and for providing the code framework, 
that allowed the numerical evaluation of our findings. The first author expresses his sincere gratitude to his friend
and colleague Max Scharf, as well as to Andr\'{e} Uschmajew and Roland Herzog for their valuable feedback and insightful discussions.

\newpage

\bibliographystyle{siamplain}
\bibliography{TTN}
\end{document}